\documentclass[10pt, a4paper]{article}
\usepackage[T1]{fontenc}
\usepackage{times}
\usepackage{enumitem}
\usepackage{amsmath, amssymb, amsthm, amsfonts}
\usepackage{graphicx}
\usepackage{caption, subcaption}
\theoremstyle{underline}
\newtheorem{definition}{Definition}
\newtheorem{theorem}{Theorem}

\newtheorem{lemma}{Lemma}

\newtheorem{remark}{Remark}

\usepackage[margin=2cm]{geometry}
\usepackage[onehalfspacing]{setspace}

\newcounter{tmp}
\usepackage{setspace}
\usepackage[numbers]{natbib}
\usepackage{hyperref}

\title{\vspace{-15mm}   On The Relationship Between The Logarithmic Lower
Order of Coefficients and The Growth of Solutions of Complex Linear
Differential Equations in $\overline{\mathbb{C}}\setminus\{z_{0}\}$ }
\author{Abdelkader Dahmani and Benharrat Bela\"{\i}di\footnote{%
Corresponding author }\\ \\Department of Mathematics, Laboratory of Pure and Applied Mathematics,\\
	University of Mostaganem (UMAB), B. P. 227 Mostaganem-Algeria\\
	abdelkader.dahmani.etu@univ-mosta.dz
	\\benharrat.belaidi@univ-mosta.dz}
\date{}
\begin{document}
\maketitle

\begin{abstract}
In this article, we study  the growth of  solutions of the homogeneous   complex linear
differential equation
\begin{equation*}
	f^{(k)}+A_{k-1}(z)f^{(k-1)}+\cdots+A_{1}(z)f^{\prime}+ A_{0}(z)f=0,
\end{equation*}%
where the coefficients $A_{j}(z)$ $(j=0,1,\ldots ,k-1)$ are analytic or meromorphic functions in $\overline{\mathbb{C}}\setminus\{z_{0}\}$. Under the sufficient condition that there exists one dominant coefficient by its logarithmic lower order or by its logarithmic lower type. We extend some precedent results due to Liu, Long and Zeng and others.
\\
	
\end{abstract}


\section{Introduction and Main Results}
In the year 2016, Fettouch and Hamouda in \cite{F2} firstly investigated the growth of solutions of linear differential equations where the coefficients are analytic or meromorphic functions in the extended complex plane except a finite singular point $\overline{\mathbb{C}}\setminus\{z_{0}\}.$ After that, and by adopting the same idea,  many valuable results were achieved on the growth of solutions of  linear differential equations (see e.g.  \cite{C2, C3, D1, D2, D3, F3, H1, L2, L3}).  The concepts of logarithmic order and logarithmic type of  Chern, \cite{C4, C5} have always been  used to study the growth of solutions to linear differential equations and linear difference equations  for that  case when the coefficients are entire or meromorphic in $\mathbb{C}$   of zero order (see e.g. \cite{B2, B3, B4, C1, F1}). Among these studies Ferraoun and Bela\"{\i}di in \cite{F1} compared the logarithmic lower order of coefficients with  the growth of solutions of higher order linear differential equations. Inspired by this work, in this article we introduce the logarithmic lower order and the logarithmic lower type of meromorphic function in   $\overline{\mathbb{C}}\setminus\{z_{0}\}$ and use them to investigate the growth of solutions to the same higher order linear differential equation and its relation with the growth of coefficients which are analytic or meromorphic in $\overline{\mathbb{C}}\setminus\{z_{0}\}$. Here we  also use  some fundamental results and some standard notations of the Nevanlinna distribution theory of meromorphic functions (see \cite{H2, L1, Y1}).  At first, let's recall some notations and definitions which can be found in \cite{D1,F2, L3}, for all $R\in (0,\infty)$ and $p\geq1,$  we define   $\exp_1R=e^R, \exp_{p+1}R=\exp(\exp_pR),$ $\log_1R=\log R$ and $ \log_{p+1}R=\log(\log_pR).$ Let $f$ be a meromorphic function in  $\overline{\mathbb{C}}-\{z_{0}\}$ and let $n(t,f)$ be the number of poles of  $f$ in $\{z\in \mathbb{C}:t\leq
 |z-z_{0}|\}\cup \{\infty \}$ ( counting multiplicities), where  $\overline{n}(t,f)$ is the number of distinct poles of  $f$. Subsequently,  the counting function of $f$ near $z_{0}$ is defined by
 \begin{equation*}
 	N_{z_{0}}(r,f)=-\int_{\infty }^{r}\frac{n(t,f)-n(\infty ,f)}{t}dt-n(\infty
 	,f)\log r
 \end{equation*}%
and
\begin{equation*}
	\overline{N}_{z_{0}}(r,f)=-\int_{\infty }^{r}\frac{\overline{n}(t,f)-\overline{n}(\infty ,f)}{t}dt-\overline{n}(\infty
	,f)\log r,
\end{equation*}%
where the
proximity function of $f$ near $z_{0}$ is defined by
\begin{equation*}
	m_{z_{0}}(r,f)=\frac{1}{2\pi }\int_{o}^{2\pi }\log ^{+}|f(z_{0}-re^{i\phi
	})|d\phi .
\end{equation*}
The characteristic function of $f$ near $z_{0}$ is defined by
\begin{equation*}
	T_{z_{0}}(r,f)=m_{z_{0}}(r,f)+N_{z_{0}}(r,f).
\end{equation*}
\hypertarget{definition 1.1}{
\begin{definition}[ \cite{L3}]
	Let $f$ be a meromorphic function in $\overline{\mathbb{%
			C}}-\{z_{0}\}$, $p$ and $q$ be two integers with $p\geq q\geq 1$. The  $%
	[p,q]$-order and the $[p,q]$-type of $f$ near $z_{0}$ are respectively  defined by%
	\begin{equation*}
		\sigma _{[ p,q]}(f,z_{0})=\limsup_{r\longrightarrow 0}\frac{\log
			_{p}^{+}T_{z_{0}}(r,f)}{\log _{q}\frac{1}{r}},\quad\tau _{[p,q]}(f,z_{0})=\limsup_{r\longrightarrow 0}\frac{\log
			_{p-1}^{+}T_{z_{0}}(r,f)}{(\log _{q-1}\frac{1}{r})^{\sigma _{[ p,q]}(f,z_{0}) }}\quad if \quad\sigma _{[ p,q]}(f,z_{0})\in (0,+\infty).
	\end{equation*}
If $f$ is an analytic  function in $\overline{\mathbb{%
		C}}-\{z_{0}\},$ then
	\begin{equation*}
	\sigma _{[ p,q]}(f,z_{0})=\limsup_{r\longrightarrow 0}\frac{\log _{p+1}^{+}M_{z_{0}}(r,f)}{\log _{q}%
			\frac{1}{r}},\quad
		\tau _{[p,q],M}(f,z_{0})=\limsup_{r\longrightarrow 0}\frac{\log _{p}^{+}M_{z_{0}}(r,f)}{(\log _{q-1}%
		\frac{1}{r})^{\sigma _{[ p,q]}(f,z_{0}) }}
	\end{equation*}%
if $\sigma _{[ p,q]}(f,z_{0})\in (0,+\infty)$, where $M_{z_{0}}(r,f)=\max \{|f(z)|:|z-z_{0}|=r\}$.	
\end{definition}}
	\hypertarget{definition 1.2}{
		\begin{definition}[\cite{D1}]
			Let $f$ be a meromorphic function in $\overline{\mathbb{%
					C}}-\{z_{0}\}$, $p$ and $q$ be two integers with $p\geq q\geq 1$. The  lower $%
			[p,q]$-order and the lower $[p,q]$-type of $f$ near $z_{0}$ are respectively  defined by%
			\begin{equation*}
				\mu _{[ p,q]}(f,z_{0})=\liminf_{r\longrightarrow 0}\frac{\log
					_{p}^{+}T_{z_{0}}(r,f)}{\log _{q}\frac{1}{r}},\quad\underline{\tau} _{[p,q]}(f,z_{0})=\liminf_{r\longrightarrow 0}\frac{\log
					_{p-1}^{+}T_{z_{0}}(r,f)}{(\log _{q-1}\frac{1}{r})^{\mu _{[ p,q]}(f,z_{0}) }}
			\end{equation*}
            if $\mu_{[ p,q]}(f,z_{0})\in (0,+\infty)$. If $f$ is an analytic  function in $\overline{\mathbb{%
					C}}-\{z_{0}\}$, then
			\begin{equation*}
				\mu _{[ p,q]}(f,z_{0})	=\liminf_{r\longrightarrow 0}\frac{\log _{p+1}^{+}M_{z_{0}}(r,f)}{\log _{q}%
					\frac{1}{r}},\quad
				\underline{\tau}_{[p,q],M}(f,z_{0})
				=\liminf_{r\longrightarrow 0}\frac{\log _{p}^{+}M_{z_{0}}(r,f)}{(\log _{q-1}%
					\frac{1}{r})^{\mu_{[ p,q]}(f,z_{0}) }}\quad if \quad\mu_{[ p,q]}(f,z_{0})\in (0,+\infty).
			\end{equation*}%
	\end{definition}}
Now we introduce the  definitions of the logarithmic order, the logarithmic type,   the logarithmic lower order and the logarithmic lower type  of meromorphic  function $f$ in $\overline{\mathbb{C}}-\{z_{0}\}$ as follow.
\hypertarget{definition 1.3 }{\begin{definition}[\cite{D2}]
		Let $f$ be a meromorphic function in $\overline{\mathbb{C}}-\{z_{0}\}$,  the logarithmic  order of $f$ near $z_{0}$ is given by
		\begin{equation*}
			\sigma _{[1,2]}(f,z_{0})=\sigma _{\log}(f,z_{0})=\limsup_{r\longrightarrow 0}\frac{\log
				^{+}T_{z_{0}}(r,f)}{\log\log \frac{1}{r}}
        \end{equation*}
        and if $\sigma _{\log}(f,z_{0})\in [1,+\infty),$ then the logarithmic type of $f$ near $z_{0}$ is defined by
        \begin{equation*}
        \tau _{[1,2]}(f,z_{0})=\tau _{\log}(f,z_{0})=\limsup_{r\longrightarrow 0}\frac{T_{z_{0}}(r,f)}{(\log \frac{1}{r})^{\sigma _{\log}(f,z_{0}) }}.
		\end{equation*}%
		If $f$ is an analytic  function in $\overline{\mathbb{C}}-\{z_{0}\}$, then
		\begin{equation*}
			\sigma _{[1,2]}(f,z_{0})=\sigma _{\log}(f,z_{0})=\limsup_{r\longrightarrow 0}\frac{\log ^{+}\log ^{+}M_{z_{0}}(r,f)}{\log\log %
				\frac{1}{r}},
       \end{equation*}%
       and if $\sigma _{\log}(f,z_{0})\in [1,+\infty),$ then the logarithmic type of $f$ near $z_{0}$ is defined by
       \begin{equation*}
        \tau _{[1,2],M}(f,z_{0})=	\tau _{\log,M}(f,z_{0})
			=\limsup_{r\longrightarrow 0}\frac{\log ^{+}M_{z_{0}}(r,f)}{(\log %
				\frac{1}{r})^{\sigma _{\log}(f,z_{0})}}.
		\end{equation*}%
\end{definition}}
\hypertarget{definition 1.4 }{\begin{definition}
		Let $f$ be a meromorphic function in $\overline{\mathbb{C}}-\{z_{0}\}$,  the logarithmic  lower order of $f$ near $z_{0}$ is given by
		\begin{equation*}
			\mu_{[1,2]}(f,z_{0})=\mu _{\log}(f,z_{0})=\liminf_{r\longrightarrow 0}\frac{\log
				^{+}T_{z_{0}}(r,f)}{\log\log \frac{1}{r}},
       \end{equation*}%
       and if $\mu _{\log}(f,z_{0})\in [1,+\infty)$, then the logarithmic lower type of $f$ near $z_{0}$ is given by
		\begin{equation*}
       \underline{\tau} _{[1,2]}(f,z_{0})=\underline{\tau}_{\log}(f,z_{0})=\liminf_{r\longrightarrow 0}\frac{T_{z_{0}}(r,f)}{(\log \frac{1}{r})^{\mu _{\log}(f,z_{0})}}.
		\end{equation*}%
		If $f$ is an analytic  function in $\overline{\mathbb{C}}-\{z_{0}\}$, then
		\begin{equation*}
			\mu _{[1,2]}(f,z_{0})=\mu _{\log}(f,z_{0})=\liminf_{r\longrightarrow 0}\frac{\log ^{+}\log ^{+}M_{z_{0}}(r,f)}{\log\log %
				\frac{1}{r}},
		\end{equation*}%
       and if $\mu _{\log}(f,z_{0})\in [1,+\infty)$, then the logarithmic lower type of $f$ near $z_{0}$ is given by
       \begin{equation*}
       \underline{\tau} _{[1,2],M}(f,z_{0})=	\underline{\tau} _{\log,M}(f,z_{0})
			=\liminf_{r\longrightarrow 0}\frac{\log ^{+}M_{z_{0}}(r,f)}{(\log %
				\frac{1}{r})^{\mu _{\log}(f,z_{0}) }}.
		\end{equation*}%
\begin{remark}
	According to (Lemma 2.2, \cite{F2}) if $f$ is a non constant meromorphic function in  $\overline{\mathbb{C}}-\{z_{0}\}$, then $g(\omega)=f(z_0-\frac{1}{\omega})$ is meromorphic in $\mathbb{C}$ and they  satisfy
	\begin{equation*}
		T(R, g)=T_{z_0}(\frac{1}{R},f).	
	\end{equation*}
	Consequently all the properties of the logarithmic order for the meromorphic functions in $\mathbb{C}$ are hold, such as all the non-constant rational functions which are analytic in  $\overline{\mathbb{C}}-\{z_{0}\}$ are of logarithmic order equals one, where there is no transcendental meromorphic function in $\overline{\mathbb{C}}-\{z_{0}\}$ of logarithmic order less than one, further, constant functions have zero logarithmic order and there are no meromorphic functions in $\overline{\mathbb{C}}-\{z_{0}\}$ of logarithmic order between zero and one (see, \cite{C4,C5}).
\end{remark}
\end{definition}}
\hypertarget{definition 1.5 }{\begin{definition}
		The  $[p,q]$ exponent of convergence of the sequence of zeros and distinct zeros of a meromorphic  function $f$ in $\overline{\mathbb{C}}-\{z_{0}\}$  are respectively  defined by
		\begin{equation*}
			\lambda _{[p,q]}(f,z_0)=\limsup_{r\longrightarrow 0}\frac{\log _{p}^{+}N_{z_{0}}(r,\frac{1}{f})}{\log _{q}%
				\frac{1}{r}},\quad\quad \overline{\lambda} _{[p,q]}(f,z_0)=\limsup_{r\longrightarrow 0}\frac{\log _{p}^{+}\overline{N}_{z_{0}}(r,\frac{1}{f})}{%
				\log _{q}\frac{1}{r}}.
		\end{equation*}
\end{definition}}

\hypertarget{definition 1.6 }{\begin{definition}
		The  logarithmic exponent of convergence of the sequence of zeros and distinct zeros of a meromorphic  function $f$ in $\overline{\mathbb{C}}-\{z_{0}\}$  are respectively  defined by
		\begin{equation*}
			\lambda _{\log}(f,z_0)=\limsup_{r\longrightarrow 0}\frac{\log^{+}N_{z_{0}}(r,\frac{1}{f})}{\log\log %
				\frac{1}{r}}-1,\quad\quad \overline{\lambda} _{\log}(f,z_0)=\limsup_{r\longrightarrow 0}\frac{\log^{+}\overline{N}_{z_{0}}(r,\frac{1}{f})}{\log\log %
				\frac{1}{r}}-1.
		\end{equation*}
\end{definition}}
\appendix
Before we start stating our results, we may mention some of the previous results obtained on the growth of solutions of linear differential equations  in the extended complex plane except a finite singular point $\overline{\mathbb{C}}\setminus\{z_{0}\}$ and we focus particularly on those involving the lower order. In \cite{L2},  Liu and his co-authors considered the following complex homogeneous second order linear differential equation
\begin{equation}\label{equation 1.1}
			f^{\prime\prime}+A(z)f^{\prime}+ B(z)f=0,
\end{equation}%
where $A(z)$ and $B(z)$ are both analytic functions in $\overline{\mathbb{C}}\setminus\{z_{0}\},$ and obtained the following theorems.
\begingroup
\setcounter{tmp}{\value{theorem}}
\setcounter{theorem}{0}
\renewcommand*{\thetheorem}{\Alph{theorem}}
 \hypertarget{theorem A}{
\begin{theorem}[\cite{L2}]
	Let  $A(z)$ and $B(z)$  be analytic functions in $\overline{\mathbb{C}}\backslash \{z_0\}$  satisfying $\mu(A, z_0)<\mu(B, z_0)<\infty.$
	 Then, every non trivial  solution $f(z)$ of (\ref{equation 1.1}) that is analytic in $\overline{\mathbb{C}}\backslash \{z_0\},$  satisfies $\sigma_{2}(f, z_0)\geq\mu(B, z_0).$
\end{theorem}}
 \hypertarget{theorem B}{
\begin{theorem}[\cite{L2}]
	Let  $A(z)$ and $B(z)$  be analytic functions in $\overline{\mathbb{C}}\backslash \{z_0\}$  satisfying $\mu(A, z_0)=\mu(B, z_0)$ with $\underline{\tau }_{M}(A, z_0)<\underline{\tau }_{M}(B, z_0).$
	Then, every non trivial  solution $f(z)$ of (\ref{equation 1.1}) that is analytic in $\overline{\mathbb{C}}\backslash \{z_0\},$  satisfies $\sigma_{2}(f, z_0)\geq\mu(B, z_0).$
\end{theorem}}
For more general case, in a recent paper \cite{D1}, we considered the following higher order linear differential equation
\begin{equation}\label{equation 1.2}
	f^{(k)}+A_{k-1}(z)f^{(k-1)}+\cdots+A_{1}(z)f^{\prime}+ A_{0}(z)f=0,
\end{equation}%
where $A_{j}(z)$ $(j=0,1,\ldots ,k-1)$
are analytic functions $\overline{\mathbb{C}}-\{z_{0}\}$ , and we obtained the following results.
\hypertarget{theorem C}{
\begin{theorem}[\cite{D1}]
		Let  $A_{0}(z),..., A_{k-1}(z)$  be analytic functions in $\overline{\mathbb{C}}\backslash \{z_0\}.$  Assume that  $$\max\big\{ \sigma_{[p,q]}(A_j, z_0) : j\neq 0 \big\} < \mu_{[p,q]}(A_0, z_0)\leq \sigma_{[p,q]}(A_0, z_0) < +\infty.$$  Then, every analytic  solution $f(z)(\not\equiv0)$ in  $\overline{\mathbb{C}}\backslash \{z_0\}$ of (\ref{equation 1.2}) satisfies $$ \mu_{[p,q]}(A_0, z_0)=\mu_{[p+1,q]}(f, z_0)\leq \sigma_{[p+1,q]}(f, z_0)=\sigma_{[p,q]}(A_0, z_0).$$
\end{theorem}}
\hypertarget{theorem D}{
	\begin{theorem}[\cite{D1}]
		Let  $A_{0}(z),..., A_{k-1}(z)$  be analytic functions in $\overline{\mathbb{C}}\backslash \{z_0\}.$ Assume that  $$\max\big\{ \sigma_{[p,q]}(A_j, z_0) : j\neq 0 \big\} \leq \mu_{[p,q]}(A_0, z_0)\leq \sigma_{[p,q]}(A_0, z_0) < +\infty$$ and $$\max\big\{ \tau_{[p,q],M}(A_j, z_0) : \sigma_{[p,q]}(A_j, z_0)= \mu_{[p,q]}(A_0, z_0)>0\big\} < \underline{\tau}_{[p,q],M}(A_0, z_0) < +\infty.$$  Then, every analytic  solution $f(z)(\not\equiv0)$ in  $\overline{\mathbb{C}}\backslash \{z_0\}$ of (\ref{equation 1.2}) satisfies $ \mu_{[p,q]}(A_0, z_0)=\mu_{[p+1,q]}(f, z_0)\leq \sigma_{[p+1,q]}(f, z_0)=\sigma_{[p,q]}(A_0, z_0).$
\end{theorem}}
\endgroup
\setcounter{theorem}{\thetmp}
As we mentioned before, in this article you may find several similarities between the results obtained in \cite{F1} for complex plane $\mathbb{C}$ and the present results which we also consider them as extensions to the above theorems and other results to the case that all the coefficients are zero order analytic or meromorphic functions in
$\overline{\mathbb{C}}\backslash \{z_0\}.$ The first five theorems are for the special case that the coefficients of  (\ref{equation 1.2}) are analytic or meromorphic in $\overline{\mathbb{C}}\backslash \{z_0\}$ and $A_0(z)$ dominates them, the rest of the results are for the case when the dominant coefficient is an arbitrary $A_s(z).$ In fact, we obtain the following results.
\hypertarget{theorem
	1.1}{\begin{theorem}
	Let  $A_{0}(z),..., A_{k-1}(z)$  be analytic functions in $\overline{\mathbb{C}}\backslash \{z_0\}$ of finite logarithmic order with $\max\big\{ \sigma_{\log}(A_j, z_0) : j\neq 0 \big\} \leq \mu_{\log}(A_0, z_0)\leq \sigma_{\log}(A_0, z_0) < +\infty$   and $$\sum_{\sigma_{\log}(A_j, z_0)=\mu_{\log}(A_0, z_0)\geq 1,j\not=0}\tau_{\log}(A_j, z_0)<\underline{\tau}_{\log}(A_0, z_0)< +\infty.$$	
	 Then, every analytic  solution $f(z)(\not\equiv0)$ in  $\overline{\mathbb{C}}\backslash \{z_0\}$ of (\ref{equation 1.2}) satisfies $0\leq\mu_{\log}(A_0, z_0)-1\leq\mu_{[2,2]}(f, z_0)\leq \mu_{\log}(A_0, z_0), $ and $\mu_{\log}(A_0, z_0)=\mu_{[2,2]}(f, z_0)\leq\sigma_{[2,2]}(f, z_0)=\sigma_{\log}(A_0, z_0)=\overline{\lambda}_{[2,2]}(f-\varphi, z_0)=\lambda_{[2,2]}(f-\varphi, z_0)$ if $\mu_{\log}(A_0, z_0)>1,$  where $\varphi(z)(\not\equiv0)$ is an analytic function in $\overline{\mathbb{C}}\backslash \{z_0\}$ satisfying $\sigma_{[2,2]}(\varphi, z_0)<\mu_{\log}(A_0, z_0).$
\end{theorem}}
\hypertarget{theorem
	1.2}{\begin{theorem}
		Let  $A_{0}(z),..., A_{k-1}(z)$  be analytic functions in $\overline{\mathbb{C}}\backslash \{z_0\}$ of finite logarithmic order with $\max\big\{ \sigma_{\log}(A_j, z_0) : j\neq 0 \big\} \leq \mu_{\log}(A_0, z_0)\leq \sigma_{\log}(A_0, z_0) < +\infty$   and 	$$\limsup_{r\longrightarrow 0}\frac{\sum_{j\not=0} m_{z_0}(r, A_j)}{m_{z_0}(r, A_0)}<1.$$
		Then, every analytic  solution $f(z)(\not\equiv0)$ in  $\overline{\mathbb{C}}\backslash \{z_0\}$ of (\ref{equation 1.2}) satisfies $0\leq\mu_{\log}(A_0, z_0)-1\leq\mu_{[2,2]}(f, z_0)\leq \mu_{\log}(A_0, z_0), $ and $\mu_{\log}(A_0, z_0)=\mu_{[2,2]}(f, z_0)\leq\sigma_{[2,2]}(f, z_0)=\sigma_{\log}(A_0, z_0)=\overline{\lambda}_{[2,2]}(f-\varphi, z_0)=\lambda_{[2,2]}(f-\varphi, z_0)$ if $\mu_{\log}(A_0, z_0)>1,$ where $\varphi(z)(\not\equiv0)$ is an analytic function in $\overline{\mathbb{C}}\backslash \{z_0\}$ satisfying $\sigma_{[2,2]}(\varphi, z_0)<\mu_{\log}(A_0, z_0).$
\end{theorem}}
\hypertarget{theorem
	1.3}{\begin{theorem}
		Let  $A_{0}(z),..., A_{k-1}(z)$  be meromorphic functions in $\overline{\mathbb{C}}\backslash \{z_0\}$ of finite logarithmic order with
		\begin{equation*}
\liminf_{r\longrightarrow 0}\frac{ m_{z_0}(r, A_0)}{T_{z_0}(r, A_0)}=\delta>0,\quad \max\big\{ \sigma_{\log}(A_j, z_0) : j\neq 0 \big\} \leq \mu_{\log}(A_0, z_0)\leq \sigma_{\log}(A_0, z_0) < +\infty
		\end{equation*}
and $$\sum_{\sigma_{\log}(A_j, z_0)=\mu_{\log}(A_0, z_0)\geq 1,j\not=0}\tau_{\log}(A_j, z_0)<\delta \underline{\tau}_{\log}(A_0, z_0)< +\infty.$$	
		Then, every meromorphic  solution $f(z)(\not\equiv0)$ in  $\overline{\mathbb{C}}\backslash \{z_0\}$ of (\ref{equation 1.2}) satisfies $0\leq\mu_{\log}(A_0, z_0)-1\leq\mu_{[2,2]}(f, z_0), $ and $\mu_{\log}(A_0, z_0)\leq\mu_{[2,2]}(f, z_0)$ if $\mu_{\log}(A_0, z_0)>1.$
\end{theorem}}
\hypertarget{theorem
	1.4}{\begin{theorem}
			Let  $A_{0}(z),..., A_{k-1}(z)$  be meromorphic functions in $\overline{\mathbb{C}}\backslash \{z_0\}$ of finite logarithmic order with
       \begin{equation*}
\liminf_{r\longrightarrow 0}\frac{ m_{z_0}(r, A_0)}{T_{z_0}(r, A_0)}=\delta>0\quad and \quad \limsup_{r\longrightarrow 0}\frac{\sum_{j\not=0} m_{z_0}(r, A_j)}{m_{z_0}(r, A_0)}<1.
       \end{equation*}
		Then, every meromorphic  solution $f(z)(\not\equiv0)$ in  $\overline{\mathbb{C}}\backslash \{z_0\}$ of (\ref{equation 1.2}) satisfies $0\leq\mu_{\log}(A_0, z_0)-1\leq\mu_{[2,2]}(f, z_0)$ and $\mu_{\log}(A_0, z_0)\leq\mu_{[2,2]}(f, z_0)$ if $\mu_{\log}(A_0, z_0)>1.$
\end{theorem}}
\hypertarget{theorem
	1.5}{\begin{theorem}
		Let  $A_{0}(z),..., A_{k-1}(z)$  be meromorphic functions in $\overline{\mathbb{C}}\backslash \{z_0\}$ of finite logarithmic order with $\lambda_{\log}(\frac{1}{A_0}, z_0)+1<\mu_{\log}(A_0, z_0),$  $\max\big\{ \sigma_{\log}(A_j, z_0) : j\neq 0 \big\} \leq \mu_{\log}(A_0, z_0)\leq \sigma_{\log}(A_0, z_0) < +\infty$   and $$\sum_{\sigma_{\log}(A_j, z_0)=\mu_{\log}(A_0, z_0),j\not=0}\tau_{\log}(A_j, z_0)< \underline{\tau}_{\log}(A_0, z_0)< +\infty.$$	
		Then, every meromorphic  solution $f(z)(\not\equiv0)$ in  $\overline{\mathbb{C}}\backslash \{z_0\}$ of (\ref{equation 1.2}) satisfies  $1<\mu_{\log}(A_0, z_0)\leq\mu_{[2,2]}(f, z_0). $
\end{theorem}}
\hypertarget{theorem
	1.6}{\begin{theorem}
		Let  $A_{0}(z),..., A_{k-1}(z)$  be analytic functions in $\overline{\mathbb{C}}\backslash \{z_0\}$ of finite logarithmic order. Suppose there exists  an integer $s$ $(0\leq s \leq k-1)$  such that  $A_{s}(z)$ satisfies $\max\big\{ \sigma_{\log}(A_j, z_0) : j\neq s \big\} \leq \mu_{\log}(A_s, z_0) < +\infty$ and $$\sum_{\sigma_{\log}(A_j, z_0)=\mu_{\log}(A_s, z_0)\geq 1,j\not=s}\tau_{\log}(A_j, z_0)<\underline{\tau}_{\log}(A_s, z_0)< +\infty.$$	
		Then, every meromorphic  solution $f(z)(\not\equiv0)$ in  $\overline{\mathbb{C}}\backslash \{z_0\}$ of (\ref{equation 1.2}) satisfies $\mu_{[2,2]}(f, z_0)-1\leq \mu_{\log}(A_s, z_0)-1\leq \mu_{\log}(f, z_0).$ Further, if $\mu_{\log}(A_s, z_0)>1,$ then $\mu_{[2,2]}(f, z_0)\leq \mu_{\log}(A_s, z_0)\leq \mu_{\log}(f, z_0).$
\end{theorem}}
\begin{remark}
	As for the case when the dominant coefficient is $A_0(z)$, we can replace the condition  $$\sum_{\sigma_{\log}(A_j, z_0)=\mu_{\log}(A_s, z_0)\geq 1,j\not=s}\tau_{\log}(A_j, z_0)<\underline{\tau}_{\log}(A_s, z_0)< +\infty$$ in \hyperlink{theorem 1.6}{Theorem 6} by $\limsup_{r\longrightarrow 0}\frac{\sum_{j\not=s} m_{z_0}(r, A_j)}{m_{z_0}(r, A_s)}<1.$
\end{remark}
\hypertarget{theorem
	1.7}{\begin{theorem}
		Let  $A_{0}(z),..., A_{k-1}(z)$  be meromorphic functions in $\overline{\mathbb{C}}\backslash \{z_0\}$ of finite logarithmic order. Suppose there exists  an integer $s$ $(0\leq s \leq k-1)$ such that  $A_{s}(z)$ satisfies $\liminf_{r\longrightarrow 0}\frac{ m_{z_0}(r, A_s)}{T_{z_0}(r, A_s)}=\delta>0,$ $\max\big\{ \sigma_{\log}(A_j, z_0) : j\neq s \big\} \leq \mu_{\log}(A_s, z_0) < +\infty$ and
		$$\sum_{\sigma_{\log}(A_j, z_0)=\mu_{\log}(A_s, z_0)\geq 1,j\not=s}\tau_{\log}(A_j, z_0)<\delta\underline{\tau}_{\log}(A_s, z_0)< +\infty.$$
		Then, every meromorphic  solution $f(z)(\not\equiv0)$ in  $\overline{\mathbb{C}}\backslash \{z_0\}$ of (\ref{equation 1.2}) satisfies $0\leq \mu_{\log}(A_s, z_0)-1\leq \mu_{\log}(f, z_0)$ and $\mu_{\log}(A_s, z_0)\leq \mu_{\log}(f, z_0)$ if $\mu_{\log}(A_s, z_0)>1.$
\end{theorem}}
\begin{remark}
	We can also replace the conditions $\max\big\{ \sigma_{\log}(A_j, z_0) : j\neq s \big\} \leq \mu_{\log}(A_s, z_0)< +\infty$ and
\begin{equation*}
\sum_{\sigma_{\log}(A_j, z_0)=\mu_{\log}(A_s, z_0)\geq 1,j\not=s}\tau_{\log}(A_j, z_0)<\delta\underline{\tau}_{\log}(A_s, z_0)< +\infty
\end{equation*}
in \hyperlink{theorem 1.7}{Theorem 7}	by $\limsup_{r\longrightarrow 0}\frac{\sum_{j\not=s} m_{z_0}(r, A_j)}{m_{z_0}(r, A_s)}<1$ or we replace the condition $\liminf_{r\longrightarrow 0}\frac{ m_{z_0}(r, A_s)}{T_{z_0}(r, A_s)}=\delta>0$ by $\lambda_{\log}(\frac{1}{A_s}, z_0)+1<\mu_{\log}(A_s, z_0),$ which clearly includes  the assumption that $\mu_{\log}(A_s, z_0)>1.$
\end{remark}
\begin{remark}
	The results in \hyperlink{theorem 1.6}{Theorem 6} and \hyperlink{theorem 1.7}{Theorem 7} may be understood as an extension respectively of Theorem 6 and Theorem 2 in \cite{D3} when an arbitrary coefficient $A_s$ dominating the others coefficients by its lower logarithmic order and lower logarithmic type.
\end{remark}
\section*{Some lemmas}
We need the following lemmas  to prove our results. First, we  define the logarithmic measure of a set $E\subset \left(0,
1\right)$ by $m_l(E)=\int_{E}\frac{dr}{r}.$
\hypertarget{lemma 2.1}{
	\begin{lemma}[\cite{D1, L3}]
		Let $f(z)$ be a non-constant analytic function in  $\overline{\mathbb{C}}\backslash \{z_0\}$ and let $V_{z_0, }(r, f)$ be a central index of $f(z)$ near the singular point $z_0.$  Then
		\begin{equation*}
		\sigma _{[ p,q]}(f,z_{0})=\limsup_{r\longrightarrow 0}\frac{\log^+_p {V_{z_0}(r, f)}}{\log_q{\frac{1}{r}}},	\quad\quad\mu_{[p,q]}(f, z_0)=	\liminf_{r\longrightarrow 0}\frac{\log^+_p {V_{z_0}(r, f)}}{\log_q{\frac{1}{r}}}.
		\end{equation*}
	\end{lemma}}
\hypertarget{lemma 2.2}{
	\begin{lemma}[\cite{H1}]
		Let $f$ be a non-constant meromorphic function in  $\overline{\mathbb{C}}\backslash \{z_0\}$ . Then  there exists a set $E_1$ of $\left(0,
		1\right)$ that has finite logarithmic measure such that for all $j=0,...,k,$ we have
		\begin{equation*}
			\frac{f^{(j)}(z_{r})}{f(z_{r})}=\left( \frac{V_{z_{0}}(r,f)}{z_{0}-z_{r}}\right) ^{j}(1+o(1)),
		\end{equation*}
		as $r\rightarrow 0,$ $r\notin E_1,$ where $z_r$ is a point in the circle $|z-z_0|=r$ that satisfies  $|f(z_r)|=\max_{|z-z_0|=r}|f(z)|.$
\end{lemma}}
\hypertarget{lemma 2.3}{
	\begin{lemma}[\cite{D1, L3}]
		Let $f$ be a non-constant analytic function in  $\overline{\mathbb{C}}\backslash \{z_0\}$ with $\mu_{[p,q]}(f, z_0)\leq\sigma_{[p,q]}(f, z_0) <\infty.$  Then
			\begin{enumerate}[label=(\roman*)]
				\item there exists a set $E_2 \subset \left(0,
				1\right)$ that has infinite logarithmic measure such that for all $|z-z_0|=r\in
				E_2,$ we have
				\begin{equation*}
					\sigma_{[p,q]}(f, z_0)=\lim_{r\longrightarrow 0}\frac{\log_p {T_{z_0}(r, f)}}{\log_q{\frac{1}{r}} }=\lim_{r\longrightarrow 0}\frac{\log_{p+1} {M_{z_0}(r, f)}}{\log_q{\frac{1}{r}} }.
				\end{equation*}
			\item there exists a set $E_3 \subset \left(0,
		1\right)$ that has infinite logarithmic measure such that for all $|z-z_0|=r\in
		E_3,$ we have
		\begin{equation*}
			\mu_{[p,q]}(f, z_0)=\lim_{r\longrightarrow 0}\frac{\log_p {T_{z_0}(r, f)}}{\log_q{\frac{1}{r}} }=\lim_{r\longrightarrow 0}\frac{\log_{p+1} {M_{z_0}(r, f)}}{\log_q{\frac{1}{r}} }.
		\end{equation*}
	\end{enumerate}
\end{lemma}}
By using  similar proofs as for  \hyperlink{lemma 2.3}{Lemma 3}, we can prove the following lemma.
\hypertarget{lemma 2.4}{
	\begin{lemma}\label{lemma 2.4}
		Let $f$ be a non-constant analytic function in  $\overline{\mathbb{C}}\backslash \{z_0\}$ with  $\mu=\mu_{\log}(f,z_0)\leq\sigma_{\log}(f,z_0)=\sigma$. Then \begin{enumerate}[label=(\roman*)]
				\item there exists a set $E_4$ of $(0, 1)$ that has infinite logarithmic measure such that for all $|z-z_0|=r\in E_4$, we have
			\begin{equation*}
				\lim_{r\longrightarrow 0}\frac{\log\log {M_{z_0}(r,f)}}{\log{\log \frac{1}{r}} }=\lim_{r\longrightarrow 0}\frac{\log {T_{z_0}(r,f)}}{\log{\log \frac{1}{r}} }=\sigma.
			\end{equation*}
			\item there exists a set $E_5$ of $(0, 1)$ that has infinite logarithmic measure such that for all $|z-z_0|=r\in E_5$, we have
		\begin{equation*}
			\lim_{r\longrightarrow 0}\frac{\log\log {M_{z_0}(r,f)}}{\log{\log \frac{1}{r}} }=\lim_{r\longrightarrow 0}\frac{\log {T_{z_0}(r,f)}}{\log{\log \frac{1}{r}} }=\mu.
		\end{equation*}
	\end{enumerate}
\end{lemma}}
\hypertarget{lemma 2.5}{
	\begin{lemma}
		Let  $A_{0}(z),...,A_{k-1}(z)$  be analytic functions in $\overline{\mathbb{C}}\backslash \{z_0\}$ of finite logarithmic order. If there exists an integer $s (0\leq s \leq k-1)$ such that    $1\leq\max\big\{\mu_{\log}(A_s, z_0), \sigma_{\log}(A_j, z_0) : j\neq s \big\} \leq \alpha,$  then every analytic  solution $f(z)(\not\equiv0)$ in  $\overline{\mathbb{C}}\backslash \{z_0\}$ of (\ref{equation 1.2}) satisfies $\mu_{[2,2]}(f, z_0)\leq \alpha.$
\end{lemma}}
\begin{proof}
	Suppose that $f(z)(\not\equiv0)$ is an analytic solution of (\ref{equation 1.2}) in  $\overline{\mathbb{C}}\backslash \{z_0\}$. By (\ref{equation 1.2}), we have
	\begin{equation}\label{equation 2.1}
		\begin{aligned}
			\bigg|\frac{f^{(k)}(z)}{f(z)}\bigg|\leq\big|A_{k-1}(z)\big|\bigg|\frac{f^{(k-1)}(z)}{f(z)}\bigg|+\cdots+\big|A_{s}(z)\big|\bigg|\frac{f^{(s)}(z)}{f(z)}\bigg|+\cdots+\big|A_{1}(z)\big|\bigg|\frac{f^{\prime}(z)}{f(z)}\bigg|+\big|A_{0}(z)\big|.
		\end{aligned}
	\end{equation} By \hyperlink{lemma 2.2}{Lemma 2}, there exists a set $E_1\subset (0,1)$ that has finite logarithmic measure such that, for all $r\not\in E_1$ and $r\rightarrow 0,$ we get
	\begin{equation}\label{equation 2.2}
		\frac{f^{(j)}(z_{r})}{f(z_{r})}=\left( \frac{V_{z_{0}}(r,f)}{z_{0}-z_{r}}\right) ^{j}(1+o(1)),\quad j=0,...,k.
	\end{equation}
	Since $\max\big\{ \sigma_{\log}(A_j, z_0) : j= 0,...,k-1, j\not = s \big\} \leq \alpha  < +\infty,$ then for any given $\varepsilon>0,$ there exists $r_1\in(0,1)$ such that for all $|z-z_0|=r\in(0,r_1),$ we have
	\begin{equation}\label{equation 2.3}
		\big|A_{j}(z)\big|\leq	M_{z_0}(r, A_j)\leq\exp\bigg\{\bigg(\log\frac{1}{r}\bigg)^{\sigma_{\log}(A_j, z_0)+\varepsilon}\bigg\}\leq\exp\bigg\{\bigg(\log\frac{1}{r}\bigg)^{\alpha+\varepsilon}\bigg\},\quad j \not=s.	
	\end{equation}
	By \hyperlink{lemma 2.4}{Lemma 4}, there exists a set $E_5\subset (0,1)$ that has infinite logarithmic measure such that, for any given $\varepsilon>0$ and for all $r\in E_5$  we get
	\begin{equation}\label{equation 2.4}
		\big|A_{s}(z)\big|\leq	M_{z_0}(r, A_s)\leq\exp\bigg\{\bigg(\log\frac{1}{r}\bigg)^{\mu_{\log}(A_s, z_0)+\varepsilon}\bigg\}\leq\exp\bigg\{\bigg(\log\frac{1}{r}\bigg)^{\alpha+\varepsilon}\bigg\}.	
	\end{equation}
	Substituting (\ref{equation 2.2}) -(\ref{equation 2.4}) into (\ref{equation 2.1}), for any given $\varepsilon>0$ and for all $r\in E_5\cap (0,r_1)\setminus E_1,$ we obtain
	\begin{equation}\label{equation 2.5}
		V_{z_0}(r,f)\leq kr\exp\bigg\{\bigg(\log\frac{1}{r}\bigg)^{\alpha+\varepsilon}\bigg\}\bigg|1+o(1)\bigg|.
	\end{equation}
	It follows by (\ref{equation 2.5}) and \hyperlink{lemma 2.1}{Lemma 1} that,  $\mu_{[2,2]}(f, z_0)\leq \alpha .$
\end{proof}
\hypertarget{lemma 2.6}{
	\begin{lemma}[\cite{D2}]
		Let  $A_{0}(z),...,A_{k-1}(z)$  be analytic functions in $\overline{\mathbb{C}}\backslash \{z_0\}$ of finite logarithmic order with $\max\big\{\sigma_{\log}(A_j, z_0) : j= 0,...,k-1 \big\} \leq \beta  < +\infty.$   Then, every analytic  solution $f(z)(\not\equiv0)$ in  $\overline{\mathbb{C}}\backslash \{z_0\}$ of (\ref{equation 1.2}) satisfies $\sigma_{[2,2]}(f, z_0)\leq \beta.$
\end{lemma}}
\hypertarget{lemma 2.7}{\begin{lemma}[\cite{F2}]
		Let $f$ be a non-constant meromorphic function in  $\overline{\mathbb{C}}\backslash \{z_0\},$ let $\lambda > 0,$ $\varepsilon > 0$ be given real constants and $j\in \mathbb{N}.$ Then
		\begin{enumerate}[label=(\roman*)]
			\item there exist a set $E_6 \subset (0,1)$ of finite logarithmic measure and a constant $C>0$ that depends only on $\lambda $ and $j$ such that for all $|z-z_0|=r\in  (0,1)\backslash E_6,$ we have
			\begin{equation}\label{equation 2.10}
				\bigg|\frac{f^{(j)}(z)}{f(z)}\bigg|\leq C\bigg[\frac{1}{r^2}T_{z_0}(\lambda r, f) \log T_{z_0}(\lambda r, f)\bigg]^{j}.
			\end{equation}
			\item there exist a set $E_7 \subset [0,2\pi)$ that has a linear  measure zero and a constant $C>0$ that depends on  $\lambda $ and $j$ such that for all $\theta \in [0,2\pi)\backslash E_7,$ there exists a constant $r_0=r_0(\theta)>0$ such that (\ref{equation 2.10}) holds for all $z$ satisfying $\arg(z-z_0) \in [0,2\pi)\backslash E_7$ and $r=|z-z_0|<r_0.$
		\end{enumerate}
\end{lemma}}
\hypertarget{lemma 2.8}{
	\begin{lemma}
			Let  $A_{0}(z),..., A_{k-1}(z)$  be analytic functions in $\overline{\mathbb{C}}\backslash \{z_0\}$ of finite logarithmic order with $\max\big\{ \sigma_{\log}(A_j, z_0) : j\neq 0 \big\} < \sigma_{\log}(A_0, z_0)= \sigma < +\infty.$   Then, every analytic  solution $f(z)(\not\equiv0)$ in  $\overline{\mathbb{C}}\backslash \{z_0\}$ of (\ref{equation 1.2}) satisfies $\sigma_{\log}(A_{0}, z_0)-1\leq\sigma_{[2,2]}(f, z_0)\leq\sigma_{\log}(A_{0}, z_0)$ and $\sigma_{[2,2]}(f, z_0)=\sigma_{\log}(A_{0}, z_0)$ if $\sigma_{\log}(A_{0}, z_0)>1.$
\end{lemma}}
\begin{proof}
	By (\ref{equation 1.2}), we have	
	\begin{equation}\label{equation 2.11}
		\begin{aligned}
			\big|A_{0}(z)\big|\leq&\bigg|\frac{f^{(k)}(z)}{f(z)}\bigg|+\big|A_{k-1}(z)\big|\bigg|\frac{f^{(k-1)}(z)}{f(z)}\bigg|+\cdots+\big|A_{1}(z)\big|\bigg|\frac{f^{\prime}(z)}{f(z)}\bigg|.
		\end{aligned}
	\end{equation}
	By \hyperlink{lemma 2.7}{Lemma 7}, there exist a set $E_6\subset (0, 1)$ having finite logarithmic measure and a constant $C>0$ that depends only on $\lambda$, such for all $r\not\in E_6,$ we have
	\begin{equation}\label{equation 2.12}
		\bigg|\frac{f^{(j)}(z)}{f(z)}\bigg|\leq C\bigg[\frac{1}{r}T_{z_0}(\lambda r, f)\bigg]^{2j}, \quad\quad j=1,...,k.
	\end{equation}
Setting $\sigma_0=\max\big\{ \sigma_{\log}(A_j, z_0) : j\neq 0 \big\}<\sigma_{\log}(A_0, z_0)= \sigma < +\infty.$ Then for any given $\varepsilon, 0<2\varepsilon<\sigma-\sigma_0,$ there exists a $r_3>0$ such that for all $|z-z_0|=r\in  (0,r_3)$ , we get
	\begin{equation}\label{equation 2.13}
		|A_{j}(z)|\leq\exp\bigg\{\bigg(\log\frac{1}{r}\bigg)^{\sigma_0+\varepsilon}\bigg\},\quad\quad j=1,...,k-1.
	\end{equation}
	By  \hyperlink{lemma 2.4}{Lemma 4}, there exists a set $E_4\subset (0,1)$ of infinite logarithmic measure such that, for all $r\in E_4$ and $|A_{0}(z)|=M_{z_0}(r, A_0)$, we have
	\begin{equation}\label{equation 2.14}
		|A_{0}(z)|>\exp\bigg\{\bigg(\log\frac{1}{r}\bigg)^{\sigma-\varepsilon}\bigg\}.
	\end{equation}
	Substituting  (\ref{equation 2.12}),  (\ref{equation 2.13}) and (\ref{equation 2.14}) into (\ref{equation 2.11}), for all $r\in E_4\cap (0,r_3)\setminus E_6,$ we obtain
	\begin{equation}\label{equation 2.15}
		\exp\bigg\{\bigg(\log\frac{1}{r}\bigg)^{\sigma-\varepsilon}\bigg\}\leq kC\bigg[\frac{1}{r}T_{z_0}(\lambda r, f)\bigg]^{2k}\exp\bigg\{\bigg(\log\frac{1}{r}\bigg)^{\sigma_0+\varepsilon}\bigg\}.
	\end{equation}
	From  (\ref{equation 2.15}) and $0<2\varepsilon<\sigma-\sigma_0,$ we get
	\begin{equation}\label{equation 2.16}
		\sigma_{[2,2]}(f,z_0)\geq\sigma-1-\varepsilon\quad if \quad\sigma\geq 1\quad and \quad\sigma_{[2,2]}(f,z_0)\geq\sigma-\varepsilon\quad if \quad\sigma> 1.
	\end{equation}
	Since $0<2\varepsilon<\sigma-\sigma_0,$ then $\sigma_{[2,2]}(f,z_0)\geq\sigma-1\geq 0$ and if $\sigma>1,$ then $\sigma_{[2,2]}(f,z_0)\geq\sigma.$ Further, by (\ref{equation 2.16}) and \hyperlink{lemma 2.6}{Lemma 6}, we have  $0\leq\sigma-1 \leq\sigma_{[2,2]}(f,z_0)\leq\sigma_{\log}(A_0, z_0)$ and  $1<\sigma \leq\sigma_{[2,2]}(f,z_0)\leq\sigma_{\log}(A_0, z_0)=\sigma.$ Thus, we obtain $\sigma_{\log}(A_0, z_0)-1\leq\sigma_{[2,2]}(f,z_0)\leq\sigma_{\log}(A_0, z_0)$ and  $1<\sigma_{\log}(A_0, z_0)=\sigma_{[2,2]}(f,z_0).$ 	
\end{proof}
\hypertarget{lemma 2.9}{
	\begin{lemma}[\cite{C2}]
		Let $f$ be a non-constant meromorphic function in $\overline{\mathbb{C}}\backslash \{z_0\}$ and let $k\in\mathbb{N}.$  Then
		\begin{equation*} m_{z_0}\big(r,\frac{f^{(k)}(z)}{f(z)}\big)= O\big(\log T_{z_0}(r,f)+\log \frac{1}{r}\big),\quad\quad \text{for all}\: r\in (0, 1)\backslash E_8
			\quad  \text{ with}      \:\: m_l(E_8)<\infty. 	
		\end{equation*}	
		If $f $ is of finite order, then
		\begin{equation*} m_{z_0}\big(r,\frac{f^{(k)}(z)}{f(z)}\big)= O\big(\log \frac{1}{r}\big),\quad\quad  r\in (0, 1).
		\end{equation*}	
\end{lemma}}
\hypertarget{lemma 2.10}{
	\begin{lemma}[\cite{L3}]
		Let $f$ be a non-constant meromorphic function in $\overline{\mathbb{C}}\backslash \{z_0\}$ and let $k$ and $j$ be two integers such that $k\not=j.$  Then
		\begin{equation*} m_{z_0}\big(r,\frac{f^{(k)}(z)}{f^{(j)}(z)}\big)\leq O\big(T_{z_0}(r,f)+\log \frac{1}{r}\big),\quad\quad \text{for all}\: r\in (0, r_4]\backslash E_9
			\quad  \text{ with} \:\: m_l(E_9)<\infty. 	
		\end{equation*}	
\end{lemma}}
\hypertarget{lemma 2.11}{
	\begin{lemma}[\cite{X1}]
	Assume  $f\not\equiv0$ is a solution of (\ref{equation 1.2}), set $g=f-\varphi.$ Then $g$  satisfies
	\begin{equation}\label{equation 2.17}
		g^{(k)}+A_{k-1}g^{(k-1)}+\cdots+A_{1}g^{\prime}+ A_{0}g=-\bigg[\varphi^{(k)}+A_{k-1}\varphi^{(k-1)}+\cdots+A_{1}\varphi^{\prime}+ A_{0}\varphi\bigg].
	\end{equation}%
\end{lemma}}
\hypertarget{lemma 2.12}{
	\begin{lemma}[\cite{F2}]
		Let $f$ be a non-constant meromorphic function in $\overline{\mathbb{C}}-\{z_{0}\}$ and set $g(\omega)=f(z_0-\frac{1}{\omega})$. Then $g(\omega)$ is meromorphic in $\mathbb{C}$ and we have $$T(R,g)=T_{z_0}(\frac{1}{R},f).$$	
\end{lemma}}
\hypertarget{lemma 2.13}{
	\begin{lemma}[\cite{B1}]
		Let $f$ be a meromorphic  function in $\mathbb{C}$ with $p\geq q\geq 1.$ Then $$\sigma_{[p,q]} (f^{\prime})=\sigma_{[p,q]} (f).$$	
\end{lemma}}
\hypertarget{lemma 2.14}{
	\begin{lemma}
		Let $f$ be a non-constant meromorphic function in $\overline{\mathbb{C}}\backslash \{z_0\}$ with $p\geq q\geq 1.$  Then  	
\begin{equation*}
			\sigma_{[p,q]} (f^{(n)}, z_0)=\sigma_{[p,q]} (f, z_0),\quad\quad n\in \mathbb{N}.
		\end{equation*} 	
\end{lemma}}
\begin{proof}
	It is sufficient to prove that $\sigma_{[p,q]} (f^{\prime}, z_0)=\sigma_{[p,q]} (f, z_0).$ By \hyperlink{lemma 2.12}{Lemma 12} , $g(\omega)=f(z_0-\frac{1}{\omega})$ is meromorphic in $\mathbb{C}$ and $\sigma_{[p,q]} (g)=\sigma_{[p,q]} (f, z_0)$. By \hyperlink{lemma 2.13}{Lemma 13} we have $\sigma_{[p,q]} (g^{\prime})=\sigma_{[p,q]} (g)$ where $f^{\prime}(z)=\frac{1}{\omega^2}g^{\prime}(\omega).$ Set $h(\omega)=\frac{1}{\omega^2}g^{\prime}(\omega).$ Clearly $\sigma_{[p,q]} (h)=\sigma_{[p,q]} (g^{\prime}).$ In the other hand by \hyperlink{lemma 2.12}{Lemma 12}, we have
	$\sigma_{[p,q]} (h)=\sigma_{[p,q]} (f^{\prime}, z_0)$. So, we deduce that $\sigma_{[p,q]} (f, z_0)=\sigma_{[p,q]} (f^{\prime}, z_0).$
\end{proof}
\hypertarget{lemma 2.15}{
	\begin{lemma}[\cite{L3}]
		Let $f$ be a non-constant meromorphic function in $\overline{\mathbb{C}}-\{z_{0}\}$ .  Then
		\begin{equation*} T_{z_0}\big(r,\frac{1}{f}\big)= T_{z_0}(r,f)+O(1). 	
		\end{equation*}	
\end{lemma}}	
\hypertarget{lemma 2.16}{
	\begin{lemma}[\cite{D2}]
		Let  $F(z)\not \equiv 0,$ $A_{0}(z),..., A_{k-1}(z)$  be analytic functions in $\overline{\mathbb{C}}\backslash \{z_0\}$ and let $f$ be a non-constant analytic solution in $\overline{\mathbb{C}}\backslash \{z_0\}$ of the equation
		\begin{equation}\label{equation 2.18}
			f^{(k)}+A_{k-1}(z)f^{(k-1)}+\cdots+A_{1}(z)f^{\prime}+ A_{0}(z)f=F(z),
		\end{equation}%
		such that  	
		\begin{equation*}
			\max\big\{\sigma_{[2,2]}(F, z_0), \sigma_{[2,2]}(A_{j}, z_0) : (j=0,...,k-1)\big\} < \sigma_{[2,2]}(f, z_0).
		\end{equation*}
		Then  $\overline{\lambda}_{[2,2]}(f, z_0)=\lambda_{[2,2]}(f, z_0)=\sigma_{[2,2]}(f, z_0).$		
\end{lemma}}

\section*{Proof of the theorems}
\subsection*{Proof of Theorem 1}
\begin{proof}
Suppose that $f(z)(\not\equiv0)$ is an analytic solution of (\ref{equation 1.2}) in  $\overline{\mathbb{C}}\backslash \{z_0\}.$  By (\ref{equation 1.2}), we get	
\begin{equation}\label{equation 3.1}
	\begin{aligned}
		-A_{0}(z)=& \frac{f^{(k)}(z)}{f(z)}+A_{k-1}(z)\frac{f^{(k-1)}(z)}{f(z)}+\cdots+A_{1}(z)\frac{f^{\prime}(z)}{f(z)}.
	\end{aligned}
\end{equation}
 By (\ref{equation 3.1}) and \hyperlink{lemma 2.9}{Lemma 9}, there exists a  set $E_8\subset (0, 1),$ having finite logarithmic measure, such that for all $r\in (0, 1)\setminus E_8,$ we have
 	\begin{equation}\label{equation 3.2}
 	\begin{aligned}
 		m_{z_0}(r,A_0(z))&\leq \sum_{j=1}^{k}m_{z_0}\big(r,\frac{f^{(j)}(z)}{f(z)}\big)+\sum_{j=1}^{k-1}m_{z_0}(r,A_{j}(z))+O(1)\\&\leq O\big(\log T_{z_0}(r,f(z))+\log \frac{1}{r}\big)+\sum_{j=1}^{k-1}m_{z_0}(r,A_{j}(z)),
 	\end{aligned}
 \end{equation}
then we obtain
	\begin{equation}\label{equation 3.3}
	T_{z_0}(r,A_0(z))=	m_{z_0}(r,A_0(z))\leq O\big(\log T_{z_0}(r,f(z))+\log \frac{1}{r}\big)+\sum_{j=1}^{k-1}T_{z_0}(r,A_{j}(z)).
\end{equation}
 First, we assume $\sigma=\max\big\{ \sigma_{\log}(A_j, z_0) : j=1,...,k-1 \big\} < \mu_{\log}(A_0, z_0)= \mu .$ Then for any given $0<2\varepsilon<\mu-\sigma,$ there exists    $r_1\in(0,1)$ such that for all $|z-z_0|=r\in(0,r_1),$ we get
	\begin{equation}\label{equation 3.4}	
	T_{z_0}(r, A_0)\geq\bigg(\log\frac{1}{r}\bigg)^{\mu_{\log}(A_0, z_0)-\varepsilon}
\end{equation}
and
\begin{equation}\label{equation 3.5}
		T_{z_0}(r, A_j)\leq\bigg(\log\frac{1}{r}\bigg)^{\sigma+\varepsilon},\quad j=1,...,k-1.
		\end{equation}
Substituting  (\ref{equation 3.4}) and  (\ref{equation 3.5}) into (\ref{equation 3.3}), for the above $\varepsilon$ and for all $|z-z_0|=r\in(0,r_1)\setminus E_8,$ we obtain
\begin{equation}\label{equation 3.6}
\bigg(\log\frac{1}{r}\bigg)^{\mu-\varepsilon}\leq O\big(\log T_{z_0}(r,f)+\log \frac{1}{r}\big)+(k-1)\bigg(\log\frac{1}{r}\bigg)^{\sigma+\varepsilon}.
\end{equation}
Then, by $0<2\varepsilon<\mu-\sigma,$ we have
\begin{equation}\label{equation 3.7}
	\big(1-o(1)\big)\bigg(\log\frac{1}{r}\bigg)^{\mu-\varepsilon}\leq O\big(\log T_{z_0}(r,f)+\log \frac{1}{r}\big).
\end{equation}
This implies that,  $\mu-1-\varepsilon\leq\mu_{[2,2]}(f,z_0).$ Since $\varepsilon>0$ is arbitrary, we obtain $0\leq\mu_{\log}(A_0, z_0)-1\leq\mu_{[2,2]}(f,z_0).$ On the other hand, by \hyperlink{lemma 2.5}{Lemma 5}, we have   $\mu_{[2,2]}(f,z_0)\leq\mu_{\log}(A_0, z_0).$ Thus, we get  that every analytic solution $f(z)(\not\equiv0)$ in  $\overline{\mathbb{C}}\backslash \{z_0\}$ of (\ref{equation 1.2}) satisfies $\mu_{\log}(A_0, z_0)-1\leq\mu_{[2,2]}(f,z_0)\leq\mu_{\log}(A_0, z_0).$ Furthermore, if $\mu_{\log}(A_0, z_0)>1, $ then by (\ref{equation 3.7}) and \hyperlink{lemma 2.5}{Lemma 5}, we obtain $\mu_{[2,2]}(f,z_0)=\mu_{\log}(A_0, z_0).$ Now we prove that  $\overline{\lambda}_{[2,2]}(f-\varphi, z_0)=\lambda_{[2,2]}(f-\varphi, z_0)=\sigma_{[2,2]}(f, z_0)=\sigma_{\log}(A_{0}, z_0)>1.$ Set $g=f-\varphi.$  By \hyperlink{lemma 2.8}{Lemma 8} and since
$\varphi(z)(\not \equiv 0)$ satisfies $\sigma_{[2,2]}(\varphi,z_0)<\mu_{\log}(A_0, z_0)\leq\sigma_{\log}(A_0, z_0),$ then we have $\sigma_{[2,2]}(g, z_0)=\sigma_{[2,2]}(f, z_0)=\sigma_{\log}(A_{0}, z_0)\geq\mu_{\log}(A_{0}, z_0)>1.$  By \hyperlink{lemma 2.11}{Lemma 11} $g$ satisfies  (\ref{equation 2.17}), Set $G=\varphi^{(k)}+A_{k-1}\varphi^{(k-1)}+\cdots+A_{1}\varphi^{\prime}+ A_{0}\varphi.$ If
 $G\equiv 0,$ then by the first part of the proof (or by \hyperlink{lemma 2.8}{Lemma 8}), we get $\sigma_{[2,2]}(\varphi,z_0)\geq\mu_{[2,2]}(\varphi,z_0)=\mu_{\log}(A_0, z_0),$ which is a contradiction, thus  $G\not\equiv 0.$ Since $G\not\equiv 0,$ then by \hyperlink{lemma 2.14}{Lemma 14}, we have  $\sigma_{[2,2]}(G, z_0)\leq\sigma_{[2,2]}(\varphi,z_0)<\mu_{\log}(A_{0}, z_0)\leq\sigma_{\log}(A_{0}, z_0)=\sigma_{[2,2]}(g, z_0).$  By \hyperlink{lemma 2.16}{Lemma 16}, we obtain
 $\overline{\lambda}_{[2,2]}(g, z_0)=\lambda_{[2,2]}(g, z_0)=\sigma_{[2,2]}(g, z_0).$ Then, we deduce that $\overline{\lambda}_{[2,2]}(f-\varphi, z_0)=\lambda_{[2,2]}(f-\varphi, z_0)=\sigma_{[2,2]}(f, z_0)=\sigma_{\log}(A_{0}, z_0).$
 \item Now we assume  that $\max\big\{ \sigma_{\log}(A_j, z_0) : j\neq 0 \big\} = \mu_{\log}(A_0, z_0)= \mu$   and
 \begin{equation*}
 \tau_1=\sum_{\sigma_{\log}(A_j, z_0)=\mu_{\log}(A_0, z_0)\geq 1,j\not=0}\tau_{\log}(A_j, z_0)<\underline{\tau}_{\log}(A_0, z_0)=\underline{\tau}.
 \end{equation*}
 Then, there exists a set $J\subseteq \{1,...,k\},$ such that for $j\in J,$ we get  $\sigma_{\log}(A_j, z_0)  = \mu_{\log}(A_0, z_0)= \mu$ with $\tau_1=\sum_{j\in J}\tau_{\log}(A_j, z_0)<\underline{\tau}_{\log}(A_0, z_0)=\underline{\tau},$ where for $j\in \{1,...,k\}\setminus J,$ we have $\sigma_{\log}(A_j, z_0)  < \mu_{\log}(A_0, z_0)= \mu.$ Hence, for any given $\varepsilon$  $\big(0<\varepsilon<\frac{\underline{\tau}-\tau_{1}}{k}\big),$  there exists a constant $r_2 \in (0,1),$   such that for all $|z-z_0|=r\in (0,r_2 )$,  we have
 \begin{equation}\label{equation 3.8}
 		T_{z_0}(r, A_j)\leq\big(\tau_{\log}(A_j, z_0)+\varepsilon \big)\bigg(\log\frac{1}{r}\bigg)^{\mu_{\log}(A_0, z_0)},\quad j\in J,
 \end{equation}
\begin{equation}\label{equation 3.9}
	T_{z_0}(r, A_j)\leq\bigg(\log\frac{1}{r}\bigg)^{\sigma_0},\quad j\in \{1,...,k\}\setminus J,\quad 0<\sigma_0<\mu
\end{equation}
and
\begin{equation}\label{equation 3.10}
	T_{z_0}(r, A_0)\geq\big(\underline{\tau}-\varepsilon \big)\bigg(\log\frac{1}{r}\bigg)^{\mu_{\log}(A_0, z_0)}.
\end{equation}
  By substituting  (\ref{equation 3.8})-(\ref{equation 3.10}) into (\ref{equation 3.3}), for the above $\varepsilon$ and for all $|z-z_0|=r\in(0,r_2)\setminus E_8,$ we get
\begin{equation}\label{equation 3.11}
	\begin{aligned}
	\big(\underline{\tau}-\varepsilon \big)\bigg(\log\frac{1}{r}\bigg)^{\mu}&\leq O\big(\log T_{z_0}(r,f)+\log \frac{1}{r}\big)+\sum_{j\in J}\big(\tau_{\log}(A_j, z_0)+\varepsilon \big)\bigg(\log\frac{1}{r}\bigg)^{\mu}+\sum_{j\in \{1,...,k\}\setminus J}\bigg(\log\frac{1}{r}\bigg)^{\sigma_0}\\&\leq O\big(\log T_{z_0}(r,f)+\log \frac{1}{r}\big)+\big(\tau_1+(k-1)\varepsilon\big)\bigg(\log\frac{1}{r}\bigg)^{\mu}+(k-1)\bigg(\log\frac{1}{r}\bigg)^{\sigma_0}
		\end{aligned}
\end{equation}
and so
\begin{equation}\label{equation 3.12}
	\big(1-o(1)\big)\big(\underline{\tau}-\tau_1-k\varepsilon\big)\bigg(\log\frac{1}{r}\bigg)^{\mu}\leq O\big(\log T_{z_0}(r,f)+\log \frac{1}{r}\big).
\end{equation}
By (\ref{equation 3.12}), it follows that
\begin{equation}\label{equation 3.13}
 0 \leq\mu_{\log}(A_0, z_0)-1\leq \mu_{[2,2]}(f,z_0).
\end{equation}
  From  (\ref{equation 3.13}) and \hyperlink{lemma 2.5}{Lemma 5}, we conclude that every analytic solution $f(z)(\not\equiv0)$ in  $\overline{\mathbb{C}}\backslash \{z_0\}$ of (\ref{equation 1.2}) satisfies $\mu_{\log}(A_0, z_0)-1\leq\mu_{[2,2]}(f,z_0)\leq\mu_{\log}(A_0, z_0).$ Furthermore, if $\mu_{\log}(A_0, z_0)>1, $ then by (\ref{equation 3.12}) and \hyperlink{lemma 2.5}{Lemma 5}, we get $\mu_{[2,2]}(f,z_0)=\mu_{\log}(A_0, z_0).$  We prove that  $\overline{\lambda}_{[2,2]}(f-\varphi, z_0)=\lambda_{[2,2]}(f-\varphi, z_0)=\sigma_{[2,2]}(f, z_0)=\sigma_{\log}(A_{0}, z_0)>1,$ similarly as in the proof for the first case.
\end{proof}
\subsection*{Proof of Theorem 2}
\begin{proof}
We assume that $\limsup_{r\longrightarrow 0}\frac{\sum_{j=1}^{k-1} m_{z_0}(r, A_j)}{m_{z_0}(r, A_0)}<\beta<1.$ Then for $r\rightarrow 0,$ we have
\begin{equation}\label{equation 3.14}
	\sum_{j=1}^{k-1} m_{z_0}(r, A_j)<\beta m_{z_0}(r, A_0).
\end{equation}
Substituting  (\ref{equation 3.14})  into (\ref{equation 3.2}), for all $r\in (0,1)\setminus E_8,$  we obtain
\begin{equation}\label{equation 3.15}
(1-\beta) T_{z_0}(r, A_0)=	(1-\beta) m_{z_0}(r, A_0)\leq O\big(\log T_{z_0}(r,f)+\log \frac{1}{r}\big).
\end{equation}	
By the definition of $\mu_{\log}(A_0, z_0)= \mu ,$ for any given $\varepsilon>0$ there exists  $r_3\in(0,1)$ such that for all $|z-z_0|=r\in(0,r_3),$  (\ref{equation 3.4}) holds. Then by  substituting  (\ref{equation 3.4}) into  (\ref{equation 3.15}), for any given $\varepsilon>0$ and for all $r\in (0,r_3)\setminus E_8,$  we get
	\begin{equation}\label{equation 3.16}
	\bigg(\log\frac{1}{r}\bigg)^{\mu-\varepsilon}\leq O\big(\log T_{z_0}(r,f)+\log \frac{1}{r}\big),
	\end{equation}
	which implies that,  $\mu-1-\varepsilon\leq\mu_{[2,2]}(f,z_0).$ Since $\varepsilon>0$ is arbitrary, we obtain
	\begin{equation}\label{equation 3.17}
	 0 \leq\mu_{\log}(A_0, z_0)-1\leq \mu_{[2,2]}(f,z_0).
	 \end{equation}
	It follows by (\ref{equation 3.17}) and  \hyperlink{lemma 2.5}{Lemma 5} that  $\mu_{\log}(A_0, z_0)-1\leq\mu_{[2,2]}(f,z_0)\leq\mu_{\log}(A_0, z_0).$ Moreover, if $\mu_{\log}(A_0, z_0)>1, $ then by (\ref{equation 3.16}) and \hyperlink{lemma 2.5}{Lemma 5}, we get $\mu_{[2,2]}(f,z_0)=\mu_{\log}(A_0, z_0).$ Similarly as in the proof of \hyperlink{theorem 1.1}{Theorem 1} we prove that  $\overline{\lambda}_{[2,2]}(f-\varphi, z_0)=\lambda_{[2,2]}(f-\varphi, z_0)=\sigma_{[2,2]}(f, z_0)=\sigma_{\log}(A_{0}, z_0)>1.$
\end{proof}
\subsection*{Proof of Theorem 3}
\begin{proof}
	Suppose that $f(z)(\not\equiv0)$ is a meromorphic  solution of (\ref{equation 1.2}) in $\overline{\mathbb{C}}\backslash \{z_0\}.$ As in the proof of \hyperlink{theorem 1.1}{Theorem 1}, first, if  $\sigma=\max\big\{ \sigma_{\log}(A_j, z_0) : j=1,...,k-1 \big\} < \mu_{\log}(A_0, z_0)= \mu ,$ then  for any given $0<2\varepsilon<\mu-\sigma,$ there exists  $r_4 \in (0,1),$   such that  for all $|z-z_0|=r\in(0,r_4),$ (\ref{equation 3.5}) holds.  By the condition  $\liminf_{r\longrightarrow 0}\frac{ m_{z_0}(r, A_0)}{T_{z_0}(r, A_0)}=\delta>0$ and the definition of $\mu_{\log}(A_0, z_0)= \mu $, for the above $\varepsilon,$ there   exists   $r_5\in(0,1)$ such that for all $|z-z_0|=r\in(0,r_5),$ we have
	\begin{equation}\label{equation 3.18}	
m_{z_0}(r, A_0)\geq	\frac{\delta}{2} T_{z_0}(r, A_0)\geq	\frac{\delta}{2}\bigg(\log\frac{1}{r}\bigg)^{\mu-\frac{\varepsilon}{2}}\geq		\bigg(\log\frac{1}{r}\bigg)^{\mu-\varepsilon}.
	\end{equation}
By substituting  (\ref{equation 3.5}) and  (\ref{equation 3.18}) into (\ref{equation 3.3}), for any given $0<2\varepsilon<\mu-\sigma,$ and for all $|z-z_0|=r\in(0,r_4)\cap(0,r_5)\setminus E_8,$ we obtain
\begin{equation}\label{equation 3.19}
	\bigg(\log\frac{1}{r}\bigg)^{\mu-\varepsilon}\leq O\big(\log T_{z_0}(r,f)+\log \frac{1}{r}\big)+(k-1)\bigg(\log\frac{1}{r}\bigg)^{\sigma+\varepsilon},
\end{equation}
that is
\begin{equation}\label{equation 3.20}
	\big(	1-o(1)\big)\bigg(\log\frac{1}{r}\bigg)^{\mu-\varepsilon}\leq O\big(\log T_{z_0}(r,f)+\log \frac{1}{r}\big).
\end{equation}
It follows that,  $0\leq\mu-1-\varepsilon\leq\mu_{[2,2]}(f,z_0)$ and  $1<\mu-\varepsilon\leq\mu_{[2,2]}(f,z_0).$ Since $\varepsilon>0$ is arbitrary, we get $0\leq \mu_{\log}(A_0, z_0)-1\leq\mu_{[2,2]}(f,z_0)$ with $1<\mu_{\log}(A_0, z_0)\leq\mu_{[2,2]}(f,z_0).$  Next, if  $\max\big\{ \sigma_{\log}(A_j, z_0) : j\neq 0 \big\} = \mu_{\log}(A_0, z_0)= \mu$   and $\tau_1=\sum_{\sigma_{\log}(A_j, z_0)=\mu_{\log}(A_0, z_0)\geq 1,j\not=0}\tau_{\log}(A_j, z_0)<\delta \underline{\tau}_{\log}(A_0, z_0)=\delta \underline{\tau},$ then by the condition  $\liminf_{r\longrightarrow 0}\frac{ m_{z_0}(r, A_0)}{T_{z_0}(r, A_0)}=\delta>0$ with $\delta\leq 1$ and  the definition of $\underline{\tau}_{\log}(A_0, z_0)= \underline{\tau},$  for any given $\varepsilon>0,$ there exists a constant $r_6\in(0,1)$ such that  for all $|z-z_0|=r\in(0,r_6),$ we have
\begin{equation}\label{equation 3.21}
	\begin{aligned}
m_{z_0}(r, A_0)\geq	\big(\delta-\varepsilon\big)T_{z_0}(r, A_0)&\geq\big(\delta-\varepsilon\big)\big(\underline{\tau}-\varepsilon \big)\bigg(\log\frac{1}{r}\bigg)^{\mu}\\
&=\big(\delta\underline{\tau}-(\underline{\tau}+\delta)\varepsilon+\varepsilon^2 \big)\bigg(\log\frac{1}{r}\bigg)^{\mu}\\
&\geq\big(\delta\underline{\tau}-(\underline{\tau}+1)\varepsilon \big)\bigg(\log\frac{1}{r}\bigg)^{\mu}.
	\end{aligned}
\end{equation}
For any given $\varepsilon$  $\big(0<(\underline{\tau}+k)\varepsilon<\delta\underline{\tau}-\tau_{1}\big),$  there exists  $r_7 \in (0,1),$   such that for all $|z-z_0|=r\in (0,r_7 )$,  (\ref{equation 3.8}) and (\ref{equation 3.9}) hold. Then, by substituting  (\ref{equation 3.8}), (\ref{equation 3.9}) and (\ref{equation 3.21}) into (\ref{equation 3.3}), for the above $\varepsilon$ and for all $|z-z_0|=r\in(0,r_6)\cap(0,r_7)\setminus E_8,$ we get
\begin{equation}\label{equation 3.22}
	\begin{aligned}
	\big(\delta\underline{\tau}-(\underline{\tau}+1)\varepsilon \big)\bigg(\log\frac{1}{r}\bigg)^{\mu}&\leq O\big(\log T_{z_0}(r,f)+\log \frac{1}{r}\big)+\sum_{j\in J}\big(\tau_{\log}(A_j, z_0)+\varepsilon \big)\bigg(\log\frac{1}{r}\bigg)^{\mu}+\sum_{j\in \{1,...,k\}\setminus J}\bigg(\log\frac{1}{r}\bigg)^{\sigma_0}\\&\leq O\big(\log T_{z_0}(r,f)+\log \frac{1}{r}\big)+\big(\tau_1+(k-1)\varepsilon\big)\bigg(\log\frac{1}{r}\bigg)^{\mu}+(k-1)\bigg(\log\frac{1}{r}\bigg)^{\sigma_0}.
	\end{aligned}
\end{equation}
So
\begin{equation}\label{equation 3.23}
	\big(1-o(1)\big)\big(\delta\underline{\tau}-\tau_1-(\underline{\tau}+k)\varepsilon \big)\bigg(\log\frac{1}{r}\bigg)^{\mu}\leq O\big(\log T_{z_0}(r,f)+\log \frac{1}{r}\big).
\end{equation}
By (\ref{equation 3.23}), we obtain
$
	0 \leq\mu_{\log}(A_0, z_0)-1\leq \mu_{[2,2]}(f,z_0)
$ with $
1 <\mu_{\log}(A_0, z_0)\leq \mu_{[2,2]}(f,z_0).
$
\end{proof}	
\subsection*{Proof of Theorem 4}
\begin{proof}
Let $f(z)(\not\equiv0)$ be a meromorphic  solution of (\ref{equation 1.2}) in $\overline{\mathbb{C}}\backslash \{z_0\}.$   For any given $\varepsilon>0,$ there exists  $r_8 \in (0,1),$   such that  for all $|z-z_0|=r\in(0,r_8),$ (\ref{equation 3.14}) and (\ref{equation 3.18}) hold. Then, by substituting (\ref{equation 3.14}) and (\ref{equation 3.18}) into (\ref{equation 3.15}),  for any given $\varepsilon>0$ and for all $|z-z_0|=r\in(0,r_8)\setminus E_8,$ we have
\begin{equation}\label{equation 3.24}
	\bigg(\log\frac{1}{r}\bigg)^{\mu-\varepsilon}\leq O\big(\log T_{z_0}(r,f)+\log \frac{1}{r}\big).
\end{equation}	
It follows that $
0 \leq\mu_{\log}(A_0, z_0)-1\leq \mu_{[2,2]}(f,z_0)
$ and $
1 <\mu_{\log}(A_0, z_0)\leq \mu_{[2,2]}(f,z_0).
$	
\end{proof}	
\subsection*{Proof of Theorem 5}
\begin{proof}
	Let $f(z)(\not\equiv0)$ be a meromorphic  meromorphic  solution of (\ref{equation 1.2}) in $\overline{\mathbb{C}}\backslash \{z_0\}.$ By (\ref{equation 3.2}), for all $r\in (0, 1)\setminus E_8,$ we obtain
	\begin{equation}\label{equation 3.25}
		\begin{aligned}
	T_{z_0}(r,A_0(z))&=	m_{z_0}(r,A_0(z))	+ N_{z_0}(r,A_0(z))\\&\leq O\big(\log T_{z_0}(r,f)+\log \frac{1}{r}\big)+\sum_{j=1}^{k-1}m_{z_0}(r,A_{j})+ N_{z_0}(r,A_0(z)).
		\end{aligned}
	\end{equation}
Also as in the proof of \hyperlink{theorem 1.1}{Theorem 1}, first, if  $\sigma=\max\big\{ \sigma_{\log}(A_j, z_0) : j=1,...,k-1 \big\} < \mu_{\log}(A_0, z_0)= \mu,$ then  for any given $\varepsilon$ $ \big(0<2\varepsilon<\mu-\sigma\big),$ there exists  $r_9 \in (0,1),$   such that  for all $|z-z_0|=r\in(0,r_9),$ (\ref{equation 3.4}) and (\ref{equation 3.5}) hold.
By the definition of $\lambda_{\log}(\frac{1}{A_0}, z_0)=\lambda,$ for any given $\varepsilon$ $\big(0<2\varepsilon<\mu-\lambda-1\big),$ there exists $r_{10} \in (0,1),$   such that  for all $|z-z_0|=r\in(0,r_{10}),$ we have
\begin{equation}\label{equation 3.26}	
	N_{z_0}(r, A_0)\leq\bigg(\log\frac{1}{r}\bigg)^{\lambda_{\log}(\frac{1}{A_0}, z_0)+1+\varepsilon}.
\end{equation}
 By substituting  (\ref{equation 3.4}), (\ref{equation 3.5}) and (\ref{equation 3.26}) into (\ref{equation 3.25}), for  sufficiently small $\varepsilon$ satisfying $0<2\varepsilon<\min \big\{\mu-\sigma,  \mu-\lambda-1\big\}$ and for all $|z-z_0|=r\in(0,r_9)\cap(0,r_{10})\setminus E_8,$ we get
 \begin{equation}\label{equation 3.27}
 	\bigg(\log\frac{1}{r}\bigg)^{\mu-\varepsilon}\leq O\big(\log T_{z_0}(r,f)+\log \frac{1}{r}\big)+(k-1)\bigg(\log\frac{1}{r}\bigg)^{\sigma+\varepsilon}+\bigg(\log\frac{1}{r}\bigg)^{\lambda+1+\varepsilon}.
 \end{equation}
 It follows that
 \begin{equation}\label{equation 3.28}
 	\big(1-o(1)\big)\bigg(\log\frac{1}{r}\bigg)^{\mu-\varepsilon}\leq O\big(\log T_{z_0}(r,f)\big),
 \end{equation}
  that is  $
 1 <\mu_{\log}(A_0, z_0)\leq \mu_{[2,2]}(f,z_0).
 $
 Now, if $\max\big\{ \sigma_{\log}(A_j, z_0) : j\neq 0 \big\} = \mu_{\log}(A_0, z_0)= \mu$   and $\tau_1=\sum_{\sigma_{\log}(A_j, z_0)=\mu_{\log}(A_0, z_0),j\not=0}\tau_{\log}(A_j, z_0)<\underline{\tau}_{\log}(A_0, z_0)=\underline{\tau},$  then for any given $\varepsilon$  $\big(0<\varepsilon<\frac{\underline{\tau}-\tau_{1}}{k}\big),$  there exists a constant $r_{11} \in (0,1),$   such that for all $|z-z_0|=r\in (0,r_{11} )$,  (\ref{equation 3.8}), (\ref{equation 3.9}) and (\ref{equation 3.10}) hold. By substituting  (\ref{equation 3.8})-(\ref{equation 3.10}) and (\ref{equation 3.26}) into (\ref{equation 3.25}), for  sufficiently small $\varepsilon$ satisfying $0<\varepsilon<\min \bigg\{\frac{\underline{\tau}-\tau_{1}}{k},  \frac{\mu-\lambda-1}{2}\bigg\}$ and for all $|z-z_0|=r\in(0,r_{10})\cup(0,r_{11})\setminus E_8,$ we obtain
 \begin{equation}\label{equation 3.29}
 	\begin{aligned}
 		\big(\underline{\tau}-\varepsilon \big)\bigg(\log\frac{1}{r}\bigg)^{\mu}\leq \:&  O\big(\log T_{z_0}(r,f)+\log \frac{1}{r}\big)+\sum_{j\in J}\big(\tau_{\log}(A_j, z_0)+\varepsilon \big)\bigg(\log\frac{1}{r}\bigg)^{\mu}\\&+\sum_{j\in \{1,...,k\}\setminus J}\bigg(\log\frac{1}{r}\bigg)^{\sigma_0}+\bigg(\log\frac{1}{r}\bigg)^{\lambda+1+\varepsilon}\\\leq \:& O\big(\log T_{z_0}(r,f)+\log \frac{1}{r}\big)+\big(\tau_1+(k-1)\varepsilon\big)\bigg(\log\frac{1}{r}\bigg)^{\mu}\\&+(k-1)\bigg(\log\frac{1}{r}\bigg)^{\sigma_0}+\bigg(\log\frac{1}{r}\bigg)^{\lambda+1+\varepsilon}.
 	\end{aligned}
 \end{equation}
 So
 \begin{equation}\label{equation 3.30}
 	\big(1-o(1)\big)\big(\underline{\tau}-\tau_1-k\varepsilon\big)\bigg(\log\frac{1}{r}\bigg)^{\mu}\leq O\big(\log T_{z_0}(r,f)\big).
 \end{equation}
From (\ref{equation 3.30}), we deduce that $
1 <\mu_{\log}(A_0, z_0)\leq \mu_{[2,2]}(f,z_0).
$
\end{proof}					
\subsection*{Proof of Theorem 6}
\begin{proof}
	Let $f(z)(\not\equiv0)$ be an  analytic solution in $\overline{\mathbb{C}}\backslash \{z_0\}$  of (\ref{equation 1.2}).  We suppose that $\mu_{\log}(f, z_0)<\infty.$ By (\ref{equation 1.2}), we have
\begin{equation}\label{equation 3.31}
	\begin{aligned}
		-A_{s}(z)&	=	\frac{f^{(k)}(z)}{f^{(s)}(z)}+A_{k-1}(z)\frac{f^{(k-1)}(z)}{f^{(s)}(z)}+\cdots+A_{s+1}(z)\frac{f^{(s+1)}(z)}{f^{(s)}(z)}+A_{s-1}(z)\frac{f^{(s-1)}(z)}{f^{(s)}(z)}+\cdots+A_{0}(z)\frac{f(z)}{f^{(s)}(z)}.
	\end{aligned}
\end{equation}
By (\ref{equation 3.31}) and \hyperlink{lemma 2.10}{Lemma 10}, there exist a  constant $r_{12} \in (0,1),$ and a set $E_9\subset (0, r_{12}],$ having finite logarithmic measure, such that for all $r\in (0,r_{12}]\setminus E_9$ , we get
\begin{equation}\label{equation 3.32}
	\begin{aligned}
		m_{z_0}(r,A_s(z))\leq &\sum_{j=0, j\not=s}^{k}m_{z_0}\big(r,\frac{f^{(j)}(z)}{f^{(s)}(z)}\big)+\sum_{j=0, j\not=s}^{k-1}m_{z_0}(r,A_{j}(z))+O(1)\\ \leq &\: O\big(T_{z_0}(r,f(z))+\log \frac{1}{r}\big)+\sum_{j=0, j\not=s}^{k-1}m_{z_0}(r,A_{j}(z)).
	\end{aligned}
\end{equation}	
Then	
\begin{equation}\label{equation 3.33}
	T_{z_0}(r,A_s(z))=	m_{z_0}(r,A_s(z))\leq  O\big(T_{z_0}(r,f(z))+\log \frac{1}{r}\big)+\sum_{j=0, j\not=s}^{k-1}T_{z_0}(r,A_{j}(z)).
\end{equation}	
 First, we suppose that  $\sigma=\max\big\{ \sigma_{\log}(A_j, z_0) : j=0,...,k-1, j\not=s \big\} < \mu_{\log}(A_s, z_0)= \mu .$ Then for any given $0<2\varepsilon<\mu-\sigma,$  there exists    $r_{13}\in(0,1)$ such that for all $|z-z_0|=r\in(0,r_{13}),$ we have
	\begin{equation}\label{equation 3.34}	
		T_{z_0}(r, A_s)\geq\bigg(\log\frac{1}{r}\bigg)^{\mu_{\log}(A_s, z_0)-\varepsilon}
	\end{equation}
	and
	\begin{equation}\label{equation 3.35}
		T_{z_0}(r, A_j)\leq\bigg(\log\frac{1}{r}\bigg)^{\sigma+\varepsilon},\quad j=0,..., s-1, s+1,...,k-1.
	\end{equation}
	Substituting  (\ref{equation 3.34}) and  (\ref{equation 3.35}) into (\ref{equation 3.33}), for the above $\varepsilon$ and for all $|z-z_0|=r\in(0,r_{12}]\cap(0,r_{13})\setminus E_9,$ we obtain
	\begin{equation}\label{equation 3.36}
		\bigg(\log\frac{1}{r}\bigg)^{\mu-\varepsilon}\leq O\big( T_{z_0}(r,f)+\log \frac{1}{r}\big)+(k-1)\bigg(\log\frac{1}{r}\bigg)^{\sigma+\varepsilon}.
	\end{equation}
	Thus,
	\begin{equation}\label{equation 3.37}
		\big(1-o(1)\big)\bigg(\log\frac{1}{r}\bigg)^{\mu-\varepsilon}\leq O\big( T_{z_0}(r,f)+\log \frac{1}{r}\big),
	\end{equation}
	which implies that,  $0\leq\mu-1-\varepsilon\leq\mu_{\log}(f,z_0).$   Since $\varepsilon>0$ is arbitrary, we get $0\leq\mu_{\log}(A_s, z_0)-1\leq \mu_{\log}(f,z_0).$   On the other hand, by \hyperlink{lemma 2.5}{Lemma 5}, we have   $\mu_{[2,2]}(f,z_0)-1\leq\mu_{\log}(A_s, z_0)-1.$ Hence, $\mu_{[2,2]}(f,z_0)-1\leq\mu_{\log}(A_s, z_0)-1\leq \mu_{\log}(f,z_0).$ Furthermore, if $\mu_{\log}(A_0, z_0)>1, $ then from (\ref{equation 3.37}) and \hyperlink{lemma 2.5}{Lemma 5}, we conclude that    $\mu_{[2,2]}(f,z_0)\leq\mu_{\log}(A_s, z_0)\leq \mu_{\log}(f,z_0).$
 Now we suppose  that $\max\big\{ \sigma_{\log}(A_j, z_0) : j=0,...,k-1, j\not=s \big\} = \mu_{\log}(A_s, z_0)= \mu $   and \\$\tau_1=\sum_{\sigma_{\log}(A_j, z_0)=\mu_{\log}(A_s, z_0)\geq 1,j\not=s}\tau_{\log}(A_j, z_0)<\underline{\tau}_{\log}(A_s, z_0)=\underline{\tau}.$ Then, there exists a set $J\subseteq \{0,1,...,k\}\setminus\{s\},$ such that for $j\in J,$ we have  $\sigma_{\log}(A_j, z_0)  = \mu_{\log}(A_s, z_0)= \mu$ with $\tau_1=\sum_{j\in J}\tau_{\log}(A_j, z_0)<\underline{\tau}_{\log}(A_s, z_0)=\underline{\tau}$ and for $j\in \{0,1,...,s-1,s+1,...,k\}\setminus J,$ we get $\sigma_{\log}(A_j, z_0)  < \mu_{\log}(A_s, z_0)= \mu.$ Hence, for any given $\varepsilon$  $\big(0<\varepsilon<\frac{\underline{\tau}-\tau_{1}}{k}\big),$  there exists $r_{14} \in (0,1),$   such that for all $|z-z_0|=r\in (0,r_{14} )$,  we have
\begin{equation}\label{equation 3.38}
	T_{z_0}(r, A_j)\leq\big(\tau_{\log}(A_j, z_0)+\varepsilon \big)\bigg(\log\frac{1}{r}\bigg)^{\mu_{\log}(A_s, z_0)},\quad j\in J,
\end{equation}
\begin{equation}\label{equation 3.39}
	T_{z_0}(r, A_j)\leq\bigg(\log\frac{1}{r}\bigg)^{\sigma_0},\quad j\in \{0,1,...,s-1,s+1,...,k\}\setminus J,\quad 0<\sigma_0<\mu
\end{equation}
and
\begin{equation}\label{equation 3.40}
	T_{z_0}(r, A_s)\geq\big(\underline{\tau}-\varepsilon \big)\bigg(\log\frac{1}{r}\bigg)^{\mu_{\log}(A_s, z_0)}.
\end{equation}
By substituting  (\ref{equation 3.38})-(\ref{equation 3.40}) into (\ref{equation 3.33}), for the above $\varepsilon$ and for all $|z-z_0|=r\in(0,r_{12}]\cap(0,r_{14})\setminus E_9,$ we get
\begin{equation}\label{equation 3.41}
	\begin{aligned}
		\big(\underline{\tau}-\varepsilon \big)\bigg(\log\frac{1}{r}\bigg)^{\mu}&\leq O\big( T_{z_0}(r,f)+\log \frac{1}{r}\big)+\sum_{j\in J}\big(\tau_{\log}(A_j, z_0)+\varepsilon \big)\bigg(\log\frac{1}{r}\bigg)^{\mu}+\sum_{j\in \{0,1,...,s-1,s+1,...,k\}\setminus J}\bigg(\log\frac{1}{r}\bigg)^{\sigma_0}\\&\leq O\big( T_{z_0}(r,f)+\log \frac{1}{r}\big)+\big(\tau_1+(k-1)\varepsilon\big)\bigg(\log\frac{1}{r}\bigg)^{\mu}+(k-1)\bigg(\log\frac{1}{r}\bigg)^{\sigma_0}
	\end{aligned}
\end{equation}
and so
\begin{equation}\label{equation 3.42}
	\big(1-o(1)\big)\big(\underline{\tau}-\tau_1-k\varepsilon\big)\bigg(\log\frac{1}{r}\bigg)^{\mu}\leq O\big( T_{z_0}(r,f)+\log \frac{1}{r}\big).
\end{equation}
 It follows that
\begin{equation}\label{equation 3.43}
	0 \leq\mu_{\log}(A_s, z_0)-1\leq \mu_{\log}(f,z_0)\quad \text{with}\quad 1<\mu_{\log}(A_s, z_0)\leq \mu_{\log}(f,z_0).
\end{equation}
Then, from  (\ref{equation 3.43}) and \hyperlink{lemma 2.5}{Lemma 5}, we deduce  that every analytic solution $f(z)(\not\equiv0)$ in  $\overline{\mathbb{C}}\backslash \{z_0\}$ of (\ref{equation 1.2}) satisfies $\mu_{[2,2]}(f,z_0)-1\leq\mu_{\log}(A_s, z_0)-1\leq\mu_{\log}(f, z_0)$ and $\mu_{[2,2]}(f,z_0)\leq\mu_{\log}(A_s, z_0)\leq\mu_{\log}(f, z_0).$
\end{proof}	
\subsection*{Proof of Theorem 7}
\begin{proof}
Here we use a similar discussion as in the proof of  \hyperlink{theorem 1.6}{Theorem 6}. Let $f(z)(\not\equiv0)$ be a  meromorphic solution in $\overline{\mathbb{C}}\backslash \{z_0\}$  of (\ref{equation 1.2}) with $\mu_{\log}(f, z_0)<\infty,$ otherwise,  the result is trivial. First, if $\sigma=\max\big\{ \sigma_{\log}(A_j, z_0) : j=0,...,k-1, j\not=s \big\} < \mu_{\log}(A_s, z_0)= \mu,$ then for any given $0<2\varepsilon<\mu-\sigma,$  there exists    $r_{15}\in(0,1)$ such that for all $|z-z_0|=r\in(0,r_{15}),$ (\ref{equation 3.35}) holds.  By   $\liminf_{r\longrightarrow 0}\frac{ m_{z_0}(r, A_s)}{T_{z_0}(r, A_s)}=\delta>0$ and the definition of $\mu_{\log}(A_s, z_0)= \mu $, for the above $\varepsilon,$ there   exists    $r_{16}\in(0,1)$ such that for all $|z-z_0|=r\in(0,r_{16}),$ we have
\begin{equation}\label{equation 3.44}	
	m_{z_0}(r, A_s)\geq	\frac{\delta}{2} T_{z_0}(r, A_s)\geq	\frac{\delta}{2}\bigg(\log\frac{1}{r}\bigg)^{\mu_{\log}(A_s, z_0)-\frac{\varepsilon}{2}}\geq		\bigg(\log\frac{1}{r}\bigg)^{\mu_{\log}(A_s, z_0)-\varepsilon}.
\end{equation}
Substituting  (\ref{equation 3.35}) and  (\ref{equation 3.44}) into (\ref{equation 3.33}), for the above $\varepsilon$ and for all $|z-z_0|=r\in(0,r_{12}]\cap(0,r_{15})\cap(0,r_{16})\setminus E_9,$ we get
\begin{equation}\label{equation 3.45}
	\bigg(\log\frac{1}{r}\bigg)^{\mu-\varepsilon}\leq O\big( T_{z_0}(r,f)+\log \frac{1}{r}\big)+(k-1)\bigg(\log\frac{1}{r}\bigg)^{\sigma+\varepsilon}.
\end{equation}
Thus
\begin{equation}\label{equation 3.46}
	\big(	1-o(1)\big)\bigg(\log\frac{1}{r}\bigg)^{\mu-\varepsilon}\leq O\big( T_{z_0}(r,f)+\log \frac{1}{r}\big).
\end{equation}
It follows that,  $0\leq\mu-1-\varepsilon\leq\mu_{\log}(f,z_0).$   Since $\varepsilon>0$ is arbitrary, we obtain $0\leq\mu_{\log}(A_s, z_0)-1\leq \mu_{\log}(f,z_0).$  Now if $\max\big\{ \sigma_{\log}(A_j, z_0) : j=0,...,k-1, j\not=s \big\} = \mu_{\log}(A_s, z_0)= \mu $   and
\begin{equation*}
\tau_1=\sum_{\sigma_{\log}(A_j, z_0)=\mu_{\log}(A_s, z_0)\geq 1,j\not=s}\tau_{\log}(A_j, z_0)<\delta\underline{\tau}_{\log}(A_s, z_0)=\delta\underline{\tau},
\end{equation*}
then by   $\liminf_{r\longrightarrow 0}\frac{ m_{z_0}(r, A_s)}{T_{z_0}(r, A_s)}=\delta>0$  with $\delta\leq 1$ and the definition of $\underline{\tau}_{\log}(A_s, z_0)= \underline{\tau},$  for any given $\varepsilon>0,$ there   exists   $r_{17}\in(0,1)$ such that  for all $|z-z_0|=r\in(0,r_{17}),$ we have

\begin{equation}\label{equation 3.47}
	\begin{aligned}
		m_{z_0}(r, A_s)\geq	\big(\delta-\varepsilon\big)T_{z_0}(r, A_s)&\geq\big(\delta-\varepsilon\big)\big(\underline{\tau}-\varepsilon \big)\bigg(\log\frac{1}{r}\bigg)^{\mu_{\log}(A_s, z_0)}\\
		&\geq\big(\delta\underline{\tau}-(\underline{\tau}+1)\varepsilon \big)\bigg(\log\frac{1}{r}\bigg)^{\mu_{\log}(A_s, z_0)}.
	\end{aligned}
\end{equation}
For any given $\varepsilon$  $\big(0<(\underline{\tau}+k)\varepsilon<\delta\underline{\tau}-\tau_{1}\big),$  there exists  $r_{18} \in (0,1),$   such that for all $|z-z_0|=r\in (0,r_{18} )$,  (\ref{equation 3.38}) and (\ref{equation 3.39}) hold. By substituting  (\ref{equation 3.38}), (\ref{equation 3.39}) and (\ref{equation 3.47}) into (\ref{equation 3.33}), for the above $\varepsilon$ and for all $|z-z_0|=r\in(0,r_{12}]\cap(0,r_{17})\cap(0,r_{18})\setminus E_9,$ we obtain
\begin{equation}\label{equation 3.48}
	\begin{aligned}
		\big(\delta\underline{\tau}-(\underline{\tau}+1)\varepsilon \big)\bigg(\log\frac{1}{r}\bigg)^{\mu}\leq &\:O\big( T_{z_0}(r,f)+\log \frac{1}{r}\big)+\sum_{j\in J}\big(\tau_{\log}(A_j, z_0)+\varepsilon \big)\bigg(\log\frac{1}{r}\bigg)^{\mu}\\&+\sum_{j\in \{0,1,...,s-1,s+1,...,k\}\setminus J}\bigg(\log\frac{1}{r}\bigg)^{\sigma_0}\\\leq &\:O\big( T_{z_0}(r,f)+\log \frac{1}{r}\big)+\big(\tau_1+(k-1)\varepsilon\big)\bigg(\log\frac{1}{r}\bigg)^{\mu}+(k-1)\bigg(\log\frac{1}{r}\bigg)^{\sigma_0},
	\end{aligned}
\end{equation}
that is
\begin{equation}\label{equation 3.49}
	\big(1-o(1)\big)\big(\delta\underline{\tau}-\tau_1-(\underline{\tau}+k)\varepsilon \big)\bigg(\log\frac{1}{r}\bigg)^{\mu}\leq O\big(T_{z_0}(r,f)+\log \frac{1}{r}\big).
\end{equation}
This implies that
$0 \leq\mu_{\log}(A_s, z_0)-1\leq \mu_{\log}(f,z_0)$ and $1 <\mu_{\log}(A_s, z_0)\leq \mu_{\log}(f,z_0).$	
\end{proof}


\begin{thebibliography}{99}
\bibitem{B1} B. Bela\"{\i}di, \textit{On the [p, q]-order of meromorphic
solutions of linear differential equations}. Acta Univ. M. Belii Ser. Math.
23 (2015), 57--69.

\bibitem{B2} B. Bela\"{\i}di, \textit{Some properties of meromorphic
solutions of logarithmic order to higher order linear difference equations}.
Bul. Acad. \c{S}tiin\c{t}e Repub. Mold. Mat. 2017, no. 1(83), 15--28.

\bibitem{B3} B. Bela\"{\i}di, \textit{Study of solutions of logarithmic
order to higher order linear differential-difference equations with
coefficients having the same logarithmic order}. Univ. Iagel. Acta Math. No.
54, (2017), 15--32.

\bibitem{B4} N. Biswas, \textit{Growth of solutions of linear
differential-difference equations with coefficients having the same
logarithmic order}. Korean J. Math. 29 (2021), no. 3, 473--481.

\bibitem{C1} T. B. Cao, K. Liu and J. Wang, \textit{On the growth of
solutions of complex differential equations with entire coefficients of
finite logarithmic order}. Math. Rep. (Bucur.) 15(65) (2013), no. 3,
249--269.

\bibitem{C2} S. Cherief and S. Hamouda, \textit{Exponent of convergence of
solutions to linear differential equations near a singular point}. Grad. J.
Math. 6 (2021), no. 2, 22--30.

\bibitem{C3} S. Cherief and S. Hamouda, \textit{Growth of solutions of a
class of linear differential equations near a singular point}. Kragujevac J.
Math. 47 (2023), no. 2, 187--201.

\bibitem{C4} T. Y. P. Chern, \textit{On meromorphic functions with finite
logarithmic order}. Trans. Amer. Math. Soc. 358 (2006), no. 2, 473--489.

\bibitem{C5} T. Y. P. Chern, \textit{On the maximum modulus and the zeros of
a transcendental entire function of finite logarithmic order}. Bull. Hong
Kong Math. Soc. 2 (1999), no. 2, 271--277.

\bibitem{D1} A. Dahmani and B. Bela\"{\i}di, \textit{Growth of solutions to
complex linear differential equations in which the coefficients are analytic
functions except at a finite singular point}. Int. J. Nonlinear Anal. Appl.
14 (2023), no. 1, 473--483.

\bibitem{D2}  A. Dahmani and B. Bela\"{\i}di, On the growth of solutions of complex linear differential equations with analytic coefficients in $\overline{\mathbb{C}}-\{z_{0}\}$ of finite logarithmic order, \textit{Vestn. Udmurtsk. Univ. Mat. Mekh. Komp. Nauki}, Vol. 33, No. 3, 2023, 416-433.

\bibitem{D3}  A. Dahmani and B. Bela\"{\i}di, Growth of solutions to complex linear differential equations with analytic or meromorphic coefficients in $\overline{\mathbb{C}}-\{z_{0}\}$  of finite logarithmic order, WSEAS Transactions on Mathematics, vol. 23, 2024, 107-117. DOI:10.37394/23206.2024.23.13

\bibitem{F1} A. Ferraoun and B. Bela\"{\i}di, \textit{Growth and oscillation
of solutions to higher order linear differential equations with coefficients
of finite logarithmic order}. Sci. Stud. Res. Ser. Math. Inform. 26 (2016),
no. 2, 115--144.

\bibitem{F2} H. Fettouch and S. Hamouda, \textit{Growth of local solutions
to linear differential equations around an isolated essential singularity}.
Electron. J. Differential Equations 2016, Paper No. 226, 10 pp.

\bibitem{F3} H. Fettouch and S. Hamouda, \textit{Local growth of solutions
of linear differential equations with analytic coefficients of finite
iterated order}. Bul. Acad. \c{S}tiin\c{t}e Repub. Mold. Mat. 2021, no.
1(95)-2(96), 69--83.

\bibitem{H1} S. Hamouda, \textit{The possible orders of growth of solutions
to certain linear differential equations near a singular point}. J. Math.
Anal. Appl. 458 (2018), no. 2, 992--1008.

\bibitem{H2} W. K. Hayman, \textit{Meromorphic functions}. Oxford
Mathematical Monographs, Clarendon Press, Oxford 1964.

\bibitem{L1} I. Laine, \textit{Nevanlinna theory and complex differential
equations}. De Gruyter Studies in Mathematics, 15. Walter de Gruyter \& Co.,
Berlin, 1993.

\bibitem{L2} Y. Liu, J. Long and S. Zeng, \textit{On relationship between
lower-order of coefficients and growth of solutions of complex differential
equations near a singular point}. Chin. Quart. J. Math. 35 (2020), no. 2,
163--170.

\bibitem{L3} J. Long and S. Zeng, \textit{On [p,q]-order of growth of
solutions of complex linear differential equations near a singular point}.
Filomat 33 (2019), no. 13, 4013--4020.

\bibitem{X1} H. Y. Xu, J. Tu and X. M. Zheng, \textit{On the hyper exponent
of convergence of zeros of }$f^{\left( j\right) }-\varphi $\textit{\ of
higher order linear differential equations}. Adv. Difference Equ. 2012,
2012:114, 16 pp.

\bibitem{Y1} C. C. Yang and H. X. Yi, \textit{Uniqueness theory of
meromorphic functions}. Mathematics and its Applications, 557. Kluwer
Academic Publishers Group, Dordrecht, 2003.
\end{thebibliography}
\end{document}